\documentclass[a4paper]{article}
\usepackage{graphicx,eurosym}
\usepackage{amsmath}
\usepackage{amsfonts}
\usepackage{amsthm}
\usepackage{amssymb}
\usepackage{subfigure}
\usepackage{multirow}
\usepackage{rotating}
\usepackage{fancybox}
\usepackage{boxedminipage}
\usepackage{version}
\usepackage[colorlinks=true]{hyperref}
\newcommand{\be}{\begin{equation}}
\newcommand{\ee}{\end{equation}}
\newcommand{\ds}{\displaystyle}
\newcommand{\NN}{\mathbb{N}}

\newcommand{\ZZ}{\mathbb{Z}}

\newcommand{\FF}{\mathbb{F}}

\theoremstyle{definition}
\newtheorem{thm}{Theorem}

\newtheorem{heur}{Heuristic}

\newtheorem{dfn}{Definition}

\newtheorem{pr}{Problem}

\begin{document}
\title{Open problems in Costas arrays}
\author{Konstantinos Drakakis\footnote{The author is also affiliated with the School of Electronic, Electrical \& Mechanical Engineering, University College Dublin, Ireland.}\\UCD CASL\\\textbf{Email:} Konstantinos.Drakakis@ucd.ie}
\makeatletter{\renewcommand*{\@makefnmark}{}\footnotetext{\textbf{Address:} UCD CASL, University College Dublin, Belfield, Dublin 4, Ireland.}\makeatother}
\maketitle

\begin{abstract}
A collection of open problems in Costas arrays is presented, classified into several categories, along with the context in which they arise.
\end{abstract}

\section{Introduction}

Costas arrays are square arrays of dots/1s and blanks/0s, such that there exists exactly one dot per row and column (that is, they are permutation arrays), and such that a) no four dots not lying on a straight line form a parallelogram, b) no four dots lying on a straight line form two equidistant pairs, and c) no three dots lying on a straight line are equidistant. They appeared for the first time in 1965 in the context of SONAR detection (\cite{C1}, and later \cite{C2} as a journal publication), when J.\ P.\ Costas, disappointed by the poor performance of SONARs, used them to describe a novel frequency hopping pattern for SONAR systems with optimal auto-correlation properties (namely auto-correlation sidelobes of height at most 1). At that stage their study was entirely empirical and application-oriented. In 1984, however, after the publication of the two main construction methods for Costas arrays by S.\ Golomb \cite{G}, still the only ones of general applicability available today, they officially acquired their present name and they became an object of mathematical interest and study.

The present author, along with his collaborators, has published numerous journal publications over the past six years on several aspects of Costas arrays, namely their theory, their properties, their enumeration, and their applications, in which many questions about them were settled. However, many old questions, along with several new ones emerging in the course of the author's research, have remained open, defying all attempts to answer them. The most important and promising amongst those are collected in this work, in the hope that it will get researchers interested in Costas arrays and willing to contribute further with their efforts.

Inevitably, the selection of the material and its presentation were influenced by the author's own preferences and views. Although a conscientious effort towards impartiality was made, the broad ground rule was set that this study is intended to be centered around Costas arrays themselves, as opposed to an object or application which Costas arrays are related to. Accordingly, the selection of the various open problems presented was based on their potential advancement of our knowledge of Costas arrays themselves, be it in the direction of their theory or their applications, rather than on the advancement of some research area (e.g.\ SONAR/RADAR systems or cryptography) through the introduction (or further application) of Costas arrays in it.

An attempt has been made to group the problems presented into families, representing thematic units. Once more, it should be stressed that this grouping is, to some extent, subjective, as it will be seen that some problems would naturally fit in several groups: in such cases, the choice was based on the problem's main context and orientation.

This is not the first time a compilation of open problems in Costas arrays is attempted. In 2006, Prof.\ S.\ Rickard presented a brief list to the IMA conference \cite{R3}, while the present author published a review on Costas arrays \cite{D10}, cataloguing some open problems as well. In 2007, S.\ Golomb and G.\ Gong published another brief work on the status of Costas arrays \cite{GG}. None of these lists, however, is of length comparable to the current one. Furthermore, these works, published at least three years ago, do not cover the significant recent developments in the subject.

\section{Why Costas arrays? Some motivation}\label{motsec}

RADAR/SONAR systems detect the the distance and velocity of targets around them by periodically transmitting a waveform $W$ and listening for reflections $R$. In an ideal noiseless environment, $R$ is an exact copy of $W$, attenuated and shifted in frequency and time. The time delay indicates the distance of the target, while the frequency shift, through the Doppler effect, its velocity (assuming that the frequency content of $W$ is narrowband enough for the Doppler effect, which is multiplicative, to be well approximated by a uniform additive shift for all frequencies).

The simplest way to detect the time and frequency shifts is by applying a matched filter: $R$ is cross-correlated with shifted versions of $W$ for various time and frequency shifts, and the pair of shifts corresponding to the maximal cross-correlation are the true shifts sought:
\be (s_t^*,s_f^*)=\underset{(s_t,s_f)}{\text{argmax}}\left|\psi(s_t,s_f)\right|, \text{where } \psi(s_t,s_f)=\int_{-\infty}^{+\infty} (F_{s_f}W)(t-s_t)R(t)dt.\ee
Here, $s_t$ and $s_f$ denote the time and frequency shifts, respectively, and the operator $F_{s_f}$ maps $W$ to the signal obtained by translating every positive frequency of $W$'s spectrum by $s_f$ to the right and every negative frequency by $s_f$ to the left (so that $F_{s_f}W$ remains real).

This simple idea fails to work in practice, because all real media are incoherent: phase delay varies with frequency, hence waveforms tend to spread while traveling in the medium \cite{C2}, so that, by the time $R$ reaches the receiver, it may look so different from $W$ that coherent processing becomes inappropriate.

J.\ P.\ Costas's idea \cite{C2} was to discard the unreliable phase information, and carry out the cross-correlation based on the energy contents of $W$ and $R$ alone. Consider a waveform of the form:
\be W(t)=A\cos\left(\phi_k+2\pi \left(f_0+\frac{a_k}{n}f_1 \right)t\right),\\ t\in\left[\frac{k-1}{n}T, \frac{k}{n}T\right],\ee
where $k\in[n]$, $T$ is the time duration of the pulse, $f_0$, $f_1$ are two predetermined frequencies, $\phi_k$ are phases suitably chosen so that the phase of $W$ is continuous in $t$ ($\phi_1=0$ may be chosen), and $a:[n]\rightarrow [n]$ is a bijection. This is a frequency hopping waveform whose instantaneous frequency is
\be f(t)=f_0+\frac{a_k}{n}f_1,\ t\in\left(\frac{k-1}{n}T, \frac{k}{n}T\right), k\in[n].\ee
$W$ is completely determined by $a$, given that $f_0$, $f_1$, and $T$ are set.

Costas's idea amounts effectively to placing an energy content detector before the matched filter, thus reconstructing $a$ from $W$, and similarly for $R$: the signals fed to the matched filter can then each be abstracted as a 2D infinite sequence, representing the time-frequency plane. This sequence is full of 0s/blanks (energy is not present), except for a $n\times n$ square that corresponds to a permutation array (whose 1s/dots denote that energy is present there). The filter overlays the two 2D sequences, shifts one with respect to the other by $v$ rows vertically and $u$ columns horizontally, and counts how many pairs of dots overlap: this is the value of the cross-correlation $\Psi$ at $(u,v)$ (compare with Definition \ref{ccdefd} below). In the absence of noise, $R$ is an exact copy of $W$, only shifted in time and frequency, so the matched filter will have found the correct shift parameters when the cross-correlation becomes equal to $n$.

When noise is present, however, some of $R$'s dots may have shifted irregularly, gone altogether missing, or even new dots may have appeared: $R$ will no longer be an exact copy of $W$, and the maximal cross-correlation will no longer be $n$. One will have no alternative than to design the filter to locate the maximal cross-correlation (whose value will no longer be known a priori) and return the shift parameters corresponding to it. This maximum may no longer be unique and/or one of the (former) sidelobes may have grown taller than the main lobe. Unfortunately, both cases result in spurious target detection.

What should the form of $a$ (or the corresponding array $A$) be in order to minimize the probability of spurious detections? In the absence of noise, the cross-correlation is just a shifted form of the autocorrelation of $A$, $A$ needs to be chosen in such a way as to suppress the height of the autocorrelation sidelobes relatively to the main lobe of height $n$ as much as possible:
\[A=A^*=\underset{A}{\text{argmin}}\max_{(u,v)}\Psi_{A,A}(u,v).\]
Choosing any pair of dots in $A$, there exists a shift (their distance vector) moving one on top of the other, so sidelobes of height 1 will exist no matter what. If, however, it is stipulated that distance vectors be unique, there will be no sidelobe of height 2 or more. This, however, is precisely the Costas property (compare with Definition \ref{cosdef} below)! Autocorrelation is known as \emph{auto-ambiguity} in the SONAR/RADAR community, and waveforms with the Costas property are said to have \emph{ideal thumbtack auto-ambiguity} \cite{C2,GT}.

Why should the optimal $A$ be a permutation array, as we assumed (summarily and without any further explanation) above? Would using twice the same frequency, or using two frequencies simultaneously, not improve the autocorrelation? Costas argued on basic engineering principles that indeed it would not \cite{C2}.

Why, finally, is the full force of the Costas property needed (if indeed it is)? Would sidelobes of height 2 or more, but still much shorter than the main lobe, not be enough? Indeed they might not, but for a reason we have not mentioned so far: even if such a signal performed well under noise conditions, as discussed above, there is still the issue of multipath interference. In the real world $R$ will most likely be the sum of attenuated and shifted noisy versions of $W$ (by different attenuation and shift parameters) representing echoes from different reflection paths, which may add up constructively at the receiver: the taller the sidelobes of $W$ are, the easier it is for them to add up to a significantly tall sidelobe in $R$.

\section{Basics on Costas arrays}

Throughout this work, $[n]$ will stand for the set of the first $n$ positive integers $\{1,\ldots,n\}$, $n\in\NN$, and obvious modifications will also be used, such as $[n]-1=\{0,1,\ldots,n-1\}$ etc. It will also be convenient to consider primes to be included in prime powers as a subset (unless explicitly excluded).

\subsection{Definition of the Costas property}\label{dcpsec}

Simply put, a Costas array is a square arrangement of dots and blanks, such that there is exactly one dot per row and column, and such that all vectors between dots are distinct.

\begin{dfn}\label{cosdef}
Let $f:[n]\rightarrow [n]$ be a bijection. $f$ is said to have the \emph{Costas property} or the \emph{distinct differences property} iff the collection of vectors $\{(i-j,f(i)-f(j)): 1\leq j<i\leq n\}$, called \emph{the distance vectors}, are all distinct, or, equivalently, iff
\be\forall i,j,k\text{ such that } 1\leq i,j,i+k,j+k\leq n: f(i+k)-f(i)=f(j+k)-f(j)\Rightarrow i=j\text{ or } k=0,\ee
in which case $f$ is called a \emph{Costas permutation}. The corresponding \emph{Costas array} $A_f$ is the square $n\times n$ array where the elements at $(f(i),i),\ i\in[n]$ are equal to 1 (dots), while the remaining elements are equal to 0 (blanks):
\be A_f=[a_{ij}]=\begin{cases} 1\text{ if } i=f(j);\\ 0\text{ otherwise}\end{cases},\ j\in[n].\ee
In view of this correspondence, the terms ``Costas arrays'' and ``Costas permutations'' will be used interchangeably. $n$ is refereed to as the \emph{order} of the Costas permutation/array.
\end{dfn}

The Costas property is invariant under rotations of the array by $90^o$, horizontal and vertical flips, and flips around the diagonals, hence a Costas array of order $n>2$ gives birth to an equivalence class (EC) that contains either eight Costas arrays, or four if the array happens to be symmetric; in the latter case, we the EC is called \emph{symmetric.} Note also that both the domain and the range of $f$ can be translated, without affecting the definition of the Costas property: in particular, considering $[n]-1$ instead of $[n]$ for both the domain and the range proves occasionally to be a more suitable convention. 
 
\begin{dfn} \label{ccdefd}
Let $f,g:[n]\rightarrow [n]$ where $n\in\NN^*$, and let $u,v\in\ZZ$. The cross-correlation between $f$ and $g$ at $(u,v)$ is defined as
\be  \label{ccdef} \Psi_{f,g}(u,v)= |\{(f(i)+v,i+u):i\in[n]\}\cap \{(g(i),i):i\in[n]\}|.\ee
\end{dfn}

Informally, the cross-correlation can be visualized in the following way: first, the two Costas arrays are placed on top of each other, and then the first is translated by $v$ units vertically and $u$ horizontally. The number of pairs of overlapping dots is the value of the cross-correlation $\Psi$ at $(u,v)$.

It should perhaps be stated explicitly that in all definitions above, addition stands for ordinary addition in $\ZZ$, not modular addition. It should also be noted that Definition \ref{cosdef} can be rephrased in the following three equivalent ways as follows: a bijection $f:[n]\rightarrow [n]$ is a Costas permutation iff
\begin{itemize}
  \item the collection of vectors $\{(f(i)-f(j),i-j):\ i,j\in[n],i>j\}$ does not contain any duplications, namely iff all vectors therein, known as the \emph{distance vectors}, are distinct, which, geometrically, means that no two of them can have both the same length and the same slope; or
  \item its auto-correlation range consists of exactly three values, namely $\{n,1,0\}$ (this is based on the description in Section \ref{motsec}).
\end{itemize}

The first variant is equivalent, in turn, to the definition given in the Introduction: no four dots can form a parallelogram, or else a pair of equal distance vectors would exist, violating the Costas property (note that the cases of four dots lying on a straight line and equidistant in pairs, or of three equidistant dots on a straight line, are considered to be limit cases of a completely ``flattened'' parallelogram and are also forbidden by the definition).

\subsection{Construction algorithms}\label{conalgsec}

There exist two two main algorithms for the construction of Costas arrays, based on the algebraic theory of finite fields. They are known as the Golomb and Welch methods, and each actually comprising several submethods. They are stated below without proofs, but all details can be found in \cite{D10,DIR,G,G2,GT}. A further semi-empirical method relying on these two, known as the Rickard method, must also be mentioned alongside them, as it is currently the only additional method that has successfully led to the discovery of previously unknown Costas arrays (four so far, namely two of order 29, one of order 36, and one of order 42) \cite{R2}. They are all presented below, classified according to their applicability (this classification, also presented in \cite{DIR}, is compatible with and an improvement of the classification presented in \cite{TRD}).

\subsubsection{Generated Costas arrays}

\begin{dfn}\label{gdef}
An infinite family of Costas arrays will be characterized as ``generated'' iff its members are constructed by an algorithm whose applicability is determined by a sufficient condition involving only the order of the array.
\end{dfn}

The following is an exhaustive list of the families of generated Costas arrays currently known:

\begin{thm}[Exponential Welch construction $W_1^\text{exp}(p,\alpha,c)$]\label{w1} Let $p$ be a prime, let $\alpha$ be a primitive root of the finite field $\FF(p)$ of $p$ elements, and let $c\in[p-1]-1$ be a constant; then, the function $f:[p-1]\rightarrow [p-1]$ where $\ds f(i)=\alpha^{i-1+c}\mod p$ is a bijection with the Costas property.
\end{thm}

Note that the inverse of a $W_1^\text{exp}$-permutation, corresponding to the transpose of the corresponding $W_1^\text{exp}$ array, is not expressible by this formula for different $\alpha$ and $c$, but rather by swapping $i$ and $f(i)$. This construction is known as the Logarithmic Welch construction $W_1^\text{log}$. The two constructions together form the Welch construction $W_1$. It can be shown that the two sets of Exponential and Logarithmic Welch arrays are disjoint for $p>5$ \cite{DGO}: in particular, $W_1$-arrays for $p>5$ are never symmetric. As there are $p-1$ ways to choose $c$ and $\phi(p-1)$ ways to choose $\alpha$ ($\phi$ stands for Euler's function), there are $2(p-1)\phi(p-1)$ $W_1$-arrays in total. The property that cyclic shifts of a $W_1^\text{exp}$ permutation (or, equivalently, of the columns of a $W_1^\text{exp}$ array) lead to a new $W_1^\text{exp}$ array is known as \emph{single periodicity}. Furthermore, it can be directly verified that, if $f$ is a $W_1^\text{exp}$ permutation, then $\ds f(i)+f\left(i+\frac{p-1}{2}\right)=p$ for all $i\in\left[\frac{p-1}{2}\right]$ (or, equivalently, the right half of a $W_1^\text{exp}$ array is the same as its left half, only flipped upside down); this is known as \emph{anti-reflective symmetry}.

\begin{thm}[Exponential Welch construction $W_2^\text{exp}(p,\alpha)$]\label{w2} Let $p$ be a prime, and let $\alpha$ be a primitive root of the finite field $\FF(p)$ of $p$ elements; then, the corresponding $W_1(p,\alpha,0)$ function $g: [p-1]\to [p-1]$ satisfies $g(1)=1$ and, consequently, the function $f: [p-3]\to [p-3]$ where $f(i)=g(i+1)-1$, $i\in[p-2]$, is a bijection with the Costas property.
\end{thm}

\begin{thm}[Golomb construction $G_2(p,m,\alpha,\beta)$]\label{g2} Let $p$ be a prime, $m\in\NN$, and let $\alpha,\ \beta$ be primitive roots of the finite field $\FF(p^m)$ of $q=p^m$ elements; then, the function $f:[q-2]\rightarrow [q-2]$ where $\ds \alpha^{f(i)}+\beta^i=1$ is a bijection with the Costas property.
\end{thm}

For a given $q$, there are $\phi^2(q-1)/m$ $G_2$-arrays in total, and they are symmetric iff either $\alpha=\beta$, this special case being known as the \emph{Lempel construction} \cite{D10}, or else whenever $q=r^2$ and $\beta=\alpha^r$ \cite{DGO}.

\begin{thm}[Golomb construction $G_3(p,m,\alpha,\beta)$]\label{g3} Let $p$ be a prime, $m\in\NN$, and let $\alpha,\ \beta$ be primitive roots of the finite field $\FF(p^m)$ of $q=p^m$ elements with the property that $\ds \alpha+\beta=1$; then, the corresponding $G_2(p,m,\alpha,\beta)$ function $g:[q-2]\rightarrow [q-2]$ satisfies $g(1)=1$, and, consequently, the function $f:[q-3]\rightarrow [q-3]$ where $f(i)=g(i+1)-1$, $i\in[q-3]$, is a bijection with the Costas property.
\end{thm}

\begin{thm}[Golomb construction $G_4(m,\alpha,\beta)$]\label{g4} Let $m\in\NN$, and let $\alpha,\ \beta$ be primitive roots of the finite field $\FF(2^m)$ of $q=2^m$ elements with the property that $\ds \alpha+\beta=1$; then, the corresponding $G_2(2,m,\alpha,\beta)$ function $g:[q-2]\rightarrow [q-2]$ satisfies $g(1)=1,\ g(2)=2$, and, consequently, the function $f:[q-4]\rightarrow [q-4]$ where $f(i)=g(i+2)-2$, $i\in[q-4]$, is a bijection with the Costas property.
\end{thm}

The inclusion of $G_3$- and $G_4$-constructions amongst the generated families is far from obvious. Indeed, both rely on the existence of two primitive roots in a finite field summing up to 1, and this seems to be an applicability condition involving more parameters than the order of the finite field, in violation of Definition \ref{gdef}. It has been proved \cite{CM}, however, that every finite field contains at least one pair of primitive roots summing up to 1, so this is, in fact, not an applicability condition.

\subsubsection{Emergent Costas arrays}

\begin{dfn}
An infinite family of Costas arrays will be characterized as ``predictably emergent'' iff a) its members are constructed through a transformation of generated Costas arrays; b) the Costas property of the members of the family cannot be asserted by a condition involving the order of the array alone; and c) the Costas property of the members of the family can be asserted by a condition involving the order of the array and some additional parameters.
\end{dfn}

\begin{thm}[Exponential Welch construction $W_3^\text{exp}(p)$] \label{w3}
Let $p$ be a prime such that 2 is a primitive element of the finite field $\FF(p)$; then, the corresponding $W_1(p,2,0)$ function $g: [p-1]\to [p-1]$ satisfies $g(1)=1$, $g(2)=2$, and, consequently, the function $f: [p-3]\to [p-3]$ where $f(i)=g(i+2)-2$, $i\in[p-3]$, is a bijection with the Costas property.
\end{thm}

\begin{thm}[Golomb construction $G_4^*(p,m,\alpha,\beta)$]\label{g44} Let $p>2$ be a prime, $m\in\NN$, and let $\alpha,\ \beta$ be primitive roots of the finite field $\FF(p^m)$ of $q=p^m$ elements with the properties that $\ds \alpha+\beta=1$ and $a^2+b^{-1}=1$; then, the corresponding $G_2(p,m,\alpha,\beta)$ function $g:[q-2]\rightarrow [q-2]$ satisfies $g(1)=1,\ g(2)=q-2$, and, consequently, the function $f:[q-4]\rightarrow [q-4]$ where $f(i)=g(i+2)-1$, $i\in[q-4]$, is a bijection with the Costas property.
\end{thm}

\begin{thm}[Golomb construction $G_4^{**}(p,m,\alpha)$]\label{g444} Let $p>2$ be a prime, $m\in\NN$, and let $\alpha$ be a primitive root of the finite field $\FF(p^m)$ of $q=p^m$ elements with the property that $\ds \alpha+\alpha^2=1$; then, the corresponding $G_2(p,m,\alpha,\alpha)$ function $g:[q-2]\rightarrow [q-2]$ satisfies $g(1)=2,\ g(2)=1$, and, consequently, the function $f:[q-4]\rightarrow [q-4]$ where $f(i)=g(i+2)-2$, $i\in[q-4]$, is a bijection with the Costas property.
\end{thm}

This method is labeled as $T_4$ in \cite{G2,GT}, but a different label is used here to improve uniformity.

\begin{thm}[Golomb construction $G_5^*(p,m,\alpha,\beta)$] \label{g5} Let $p$ be a prime and $m\in\NN^*$, and let $\alpha$, $\beta$ be primitive elements of the finite field $\FF(q)$ with the properties that $\alpha+\beta=1$ and $\alpha^2+\beta^{-1}=1$; then, it can be shown that these conditions always imply that $b^2+a^{-1}=1$ as well, whence the corresponding $G_2(p,m,\alpha,\beta)$ function $g:[q-2]\rightarrow [q-2]$ satisfies $g(1)=1,\ g(2)=q-2$, and $g(q-2)=2$, and, consequently, the function $f:[q-5]\rightarrow [q-5]$ where $f(i)=g(i+3)-3$, $i\in[q-5]$, is a bijection with the Costas property.
\end{thm}

Regarding $W_3^\text{exp}$, it is not known for which primes $p$ 2 is a primitive root, although necessary conditions for this to occur can easily be found (e.g.\ 2 should not be a square). It is not even known whether 2 is a primitive root for infinitely many primes $p$ (this is a special case of Artin's Conjecture \cite{IR}). Sufficient conditions for the existence of primitive roots with the required properties for $G_4^*$, $G_4^{**}$, and $G_5^*$ are, similarly, unknown \cite{G2,GT}.

\begin{dfn}
An infinite family of Costas array will be characterized as ``unpredictably emergent'' iff a) its members are constructed through a transformation of generated Costas arrays; and b) the Costas property of the members of the family cannot be asserted by any condition other than direct verification.
\end{dfn}

\begin{heur}[Golomb construction $G_1(p,m,\alpha,\beta)$]\label{g1} Let $p$ be a prime, $m\in\NN$, $\alpha,\ \beta$ be primitive roots of the finite field $\FF(p^m)$ of $q=p^m$ elements, and let $g:[q-2]\rightarrow [q-2]$ be the corresponding $G_2(p,m,\alpha,\beta)$. It may hold true that the function $f:[q-1]\rightarrow [q-1]$ such that $f(1)=1$ and $f(i)=g(i-1)+1$, $i\in[q-2]+1$, is a bijection with the Costas property.
\end{heur}

This simply adds a corner dot to a $G_2$-Costas array.

\begin{heur}[Golomb construction $G_0(p,m,\alpha,\beta)$]\label{g0} Let $p$ be a prime, $m\in\NN$, $\alpha,\ \beta$ be primitive roots of the finite field $\FF(p^m)$ of $q=p^m$ elements, and let $g:[q-2]\rightarrow [q-2]$ be the corresponding $G_2(p,m,\alpha,\beta)$. It may hold true that the function $f:[q]\rightarrow [q]$ such that $f(1)=1$, $f(q)=q$, and $f(i)=g(i-1)+1$, $i\in[q-2]+1$, is a bijection with the Costas property.
\end{heur}

This simply adds two anti-diametrical corner dots to a $G_2$-Costas array.

\begin{heur}[Welch construction $W_0(p,\alpha,c)$]\label{w0} Let $p$ be a prime, $\alpha$ be a primitive root of the finite field $\FF(p)$ of $p$ elements, $c\in[p-1]-1$ be a constant, and let $g:[p-1]\rightarrow [p-1]$ be the corresponding $W_1(p,\alpha,c)$. It may hold true that the function $f:[p]\rightarrow [p]$ such that $f(1)=1$, and $f(i)=g(i-1)+1$, $i\in[p-1]+1$, is a bijection with the Costas property.
\end{heur}

This simply adds a corner dot to a $W_1$-Costas array.

\begin{heur}[Rickard Welch construction $RW_0(p,\alpha,c,t)$] \label{rw0}
Let $p$ be a prime, let $\alpha$ be a primitive root of the finite field $\FF(p)$ of $p$ elements, let $c\in[p-1]-1$ and $t\in[p-2]+1$ be constants, and let $g:[p-1]\rightarrow [p-1]$ be the corresponding $W_1^\text{exp}(p,\alpha,c)$ function. Consider the expansion $g: [p-1]\rightarrow [p]$, and the function $g_t:[p-1]\rightarrow [p]$ where $g_t(i)=(g(i)+t-1)\bmod p +1$. Observe that $t$ does not belong to the range of $g_t$, so it may hold true that $f:[p]\to [p]$, where $f(i)=g_t(i)$, $i\in[p-1]$, and $f(p)=t$, is a bijection with the Costas property.
\end{heur}

To understand better how this method works, start with a $(p-1)\times (p-1)$ $W_1^\text{exp}$ array and add a blank row at the bottom, thus forming a $p\times (p-1)$ array. By cyclically shifting the rows of this array $t$ times, we obtain a new $p\times (p-1)$ array where row $t$ is now blank. By appending a column to the right, then, with a dot in row $t$ we construct a new permutation array, which may have the Costas property.

\begin{heur}[Rickard Golomb construction $RG_1(p,m,\alpha,\beta, t_1,t_2)$] \label{rg1}
Let $p$ be a prime, $m\in\NN$, let $\alpha,\ \beta$ be primitive roots of the finite field $\FF(p^m)$ of $q=p^m$ elements, let $t_1,t_2\in[q-3]+1$ and let $g:[q-2]\rightarrow [q-2]$ be the corresponding $G_2(p,m,\alpha,\beta)$ function. Consider the function $g_{t_1,t_2}: [q-1]\backslash\{t_1\}\to [q-1]\backslash\{t_2\}$ where $g_{t_1,t_2}(i)=(g((i-t_1-1)\bmod (q-1)+1)+t_2-1)\bmod (q-1)+1$. It may hold true that $f:[q-1]\to [q-1]$, where $f(i)=g_{t_1,t_2}(i)$, $i\in[q-1]\backslash\{t_1\}$, and $f(t_1)=t_2$, is a bijection with the Costas property.
\end{heur}

This method can also be better understood in terms of operations on a Costas array. Start with a $(q-2)\times (q-2)$ $G_2$ array and add a blank row at the bottom and a blank column to the right, thus forming a $(q-1)\times (q-1)$ array. By cyclically shifting the columns of this array $t_1$ times and the rows $t_2$ times, we obtain a new $(q-1)\times (q-1)$ array where column $t_1$ and row $t_2$ is now blank. By adding a dot at position $(t_2,t_1)$, then, we construct a new permutation array, which may have the Costas property.

Note that Rickard Costas arrays can naturally be considered to include $G_1$ and $W_0$ Costas arrays as special cases, although in Definitions \ref{rw0} and \ref{rg1} above we excluded the parameters that would lead to such an inclusion, in order to keep the available construction methods as disjoint as possible.

\subsubsection{Sporadic Costas arrays}

\begin{dfn}
A Costas array will be characterized as ``sporadic'' iff it is neither generated nor emergent.
\end{dfn}

Sporadic Costas arrays of order $n$ exist for $6\leq n\leq 27$. A single sporadic Costas array of order 27 exists, the existence of which was first announced in \cite{DRBCIOW} and was first noted by J.K.~Beard. It is currently the largest sporadic Costas array known. It is not known whether Costas arrays of order $n$ exist for all $n\in\NN$, but, given the list of the methods above, it is clear that, unless sporadic Costas arrays exist in abundance, in particular at the orders where generated or emergent Costas arrays do not exist, this cannot happen. The smallest values of $n$ for which no Costas array of order $n$ is currently known are $n=32$ and 33 \cite{D10}.

\subsection{Known Costas arrays}\label{kcasec}

The following families of Costas arrays are known at present:
\begin{itemize}
  \item All Costas arrays of orders $n\leq 29$ have been found through exhaustive search \cite{BREMW,DIR,DIRW,DRBCIOW,RCDLW}.
  \item Costas arrays generated by the algebraic construction techniques mentioned above are available for infinitely many orders.
\end{itemize}

Table \ref{nocatbl} shows the number of known Costas arrays per order $n$. Note the main ``lobe'' formed by orders $1\leq n\leq 26$: the number of Costas arrays monotonically increases for $1\leq n\leq 16$ and monotonically decreases for $16\leq n\leq 26$. The majority of Costas arrays in the orders lying towards the interior of the lobe are sporadic Costas arrays: for example, only 16 out of the 10240 Costas arrays of order 19 are algebraically constructed (they are $W_0$-arrays, to be precise). The situation is reversed towards the edges of the lobe, where there are only two sporadic ECs for $n=26$ \cite{RCDLW} (the paper actually reports three sporadic ECs, but one of them turns out to be constructible by the Rickard method), while all Costas arrays of orders $n\leq 5$ are algebraically constructed.

\begin{table}
\centering
\begin{tabular}{|r|r||r|r||r|r||r|r|}
\hline
$n$ & $C(n)$ & $n$ & $C(n)$ & $n$ & $C(n)$ & $n$ & $C(n)$ \\\hline
1 & 1/1 & 10 & 2160/28 & 19 & 10240/12 & 28 & 712/0\\
2 & 2/1 & 11 & 4368/36 & 20 & 6464/8 & 29 & $ 164/10$\\
3 & 4/2 & 12 & 7852/34 & 21 & 3536/16 & 30 & $\geq 664/8$\\
4 & 12/2 & 13 & 12828/50 & 22 & 2052/10 & 31 & $\geq 8/0$\\
5 & 40/4 & 14 & 12752/46 & 23 & 872/20 & 32 & $\geq 0$\\
6 & 116/10 & 15 & 19612/62 & 24 & 200/0 & 33 & $\geq 0$\\
7 & 200/20 & 16 &  21104/40 & 25 & 88/4 & & \\
8 & 444/18 & 17 & 18276/38 & 26 & 56/4& & \\
9 & 760/20 & 18 & 15096/20 & 27 & 204/14& & \\
\hline
\end{tabular}
\caption{\label{nocatbl} The number of Costas arrays per order $1\leq n\leq 32$: the numbers shown after the slash are the numbers of symmetric Costas arrays (there are two symmetric Costas arrays per EC for orders $n>2$); the $\geq$ notation is used for orders not yet enumerated, to signify that there may be more arrays in addition to the ones known so far.}
\end{table}

There are two major online sources of information for Costas arrays:
\begin{itemize}
  \item Dr.\ J.K.\ Beard has prepared a database of all Costas arrays of orders $n\leq 500$, which is freely available on his website \cite{cawpb};
  \item The webpage set up by Prof.\ Scott Rickard and his associates \cite{cawp} contains not only most papers on Costas arrays published in the literature, but also a freely available Matlab toolbox on Costas arrays \cite{TRD}.
\end{itemize}

\subsection{The difference triangle}

As mentioned in Section \ref{dcpsec}, a permutation $f$ of order $n$ has the Costas property iff all distance vectors $\{(f(i)-f(j),i-j):\ i,j\in[n],i>j\}$ are distinct. In order to verify this property more easily, these vectors can be grouped together according to their second coordinate: within each such group, all first coordinate values must be distinct. In other words, setting $i-j=k$, the collection of distance vectors can be rewritten as $T(f)=\{t_k(f),\ k\in[n-1]\}$ where $t_k(f)=(f(j+k)-f(j):\ j\in[n-k])$. $T$ is said to be the \emph{difference triangle} of $f$ and $t_k$ its $k$th row, because of the way $t_k$ can be stacked one below the other to form a triangular shape, as in Table \ref{extbl}. The Costas property then is equivalent to ascertaining that no $t_k$ contains duplicate entries (this is indeed the case in Table \ref{extbl}).

\begin{table}[t]
\centering
\begin{tabular}{ccccccccccccccccccc}
1 & & 2 & & 4 & & 8 & & 5 & & 10 & & 9 & & 7 & & 3 & & 6\\\hline
& 1 & & 2 & & 4 & & -3 & & 5 & & -1 & & -2 & & -4 & & 3 & \\
& & 3 & & 6 & & 1 & & 2 & & 4 & & -3 & & -6 & & -1 & & \\
& & & 7 & & 3 & & 6 & & 1 & & 2 & & -7 & & -3 & & & \\
& & & & 4 & & 8 & & 5 & & -1 & & -2 & & -4 & & & & \\
& & & & & 9 & & 7 & & 3 & & -5 & & 1 & & & & & \\
& & & & & & 8 & & 5 & & -1 & & -2 & & & & & & \\
& & & & & & & 6 & & 1 & & 2 & & & & & & & \\
& & & & & & & & 2 & & 4 & & & & & & & & \\
& & & & & & & & & 5 & & & & & & & & &
\end{tabular}
\caption{\label{extbl} A Costas permutation along with its difference triangle}
\end{table}

According to what has been stated so far, $\ds \sum_{k=1}^{n-1}{n-k\choose 2}={n\choose 3}$ comparisons of pairs of values in the difference are required to verify the Costas property. The entries of the triangle, however, exhibit strong correlation, and this number can be drastically reduced: W.\ Chang proved in 1987 \cite{Ch} that if $t_k,\ k=1,\ldots,\ds \left\lfloor\frac{n-1}{2}\right\rfloor$ do not contain any duplication, then the remaining $t_k$ are automatically duplication-free as well. A further improvement within the first $\ds \left\lfloor\frac{n-1}{2}\right\rfloor$ was later proposed in \cite{BDR}.

\section{Core problems in the theory of Costas arrays}

\subsection{Existence of Costas arrays}

The first and foremost requirement for the well-posedness of a mathematical problem is a proof that a solution to this problem exists (or else any search for a solution would be futile). Accordingly, a requirement for the mathematical study of Costas arrays is a result determining which orders they exist in and do not exist in. In particular, do Costas arrays exist in all orders? This question was asked for the first time by S.\ Golomb and H.\ Taylor in 1984 \cite{GT} and remains open ever since:

\begin{pr}
Determine the number $C(n)$ of Costas arrays of order $n\in\NN^*$. In particular, is $C(n)>0$ for all $n>0$? If not, determine the set $X=\{n\in\NN^*:\ C(n)=0\}$.
\end{pr}

A notable attempt to compute $C(n)$ is due to W.\ Correll Jr.\ \cite{C}. Let $\mathcal{P}(n)$ stand for the set of all permutations of order $n$, let $\mathcal{D}(n)$ be the set of all pairs $(i,j)$, $1\leq i< j\leq n$, let $\mathcal{T}(n)=\{((i_1,j_1),(i_2,j_2))\in\mathcal{D}(n)\times\mathcal{D}(n):\ j_1-i_1=j_2-i_2\}$, and let $\mathcal{I}(n)=2^{\mathcal{T}(n)}$. Finally, for any $i\in \mathcal{I}(n)$, let $\mathcal{F}(i,n)=\{f\in\mathcal{P}(n):\ \forall ((i_1,j_1),(i_2,j_2))\in i,\ f(j_1)-f(i_1)=f(j_2)-f(i_2)\}$, which essentially represent ``hyperplanes'' in the space $\mathcal{P}(n)$. Note that
\be \mathcal{P}(n)\backslash\mathcal{C}(n)=\bigcup_{i\in\mathcal{I}(n)}\mathcal{F}(i,n)\leftrightarrow n!-C(n)=\left|\bigcup_{i\in\mathcal{I}(n)}\mathcal{F}(i,n)\right|,\ee
but also that elements of $\mathcal{I}(n)$ are related by inclusion:
\be \forall i_1,i_2\in\mathcal{I}(n),\ i_1\cap i_2\subset i_j\subset i_1\cup i_2,\ j=1,2.\ee
Such a structure is known as a (finite) lattice \cite{S}, and it allows one to express $|\cup_{i\in\mathcal{I}(n)}\mathcal{F}(i,n)|$ in terms of the individual $|\mathcal{F}(i,n)|$ through the M\"obius Inversion Formula \cite{S}, which is simply a generalized version of the well known Inclusion-Exclusion Principle. Unfortunately, neither the exact lattice structure of $\mathcal{I}(n)$ nor the exact value of $|\mathcal{F}(i,n)|$ for an arbitrary $i\in\mathcal{I}(n)$ are simple to compute, so this approach remains mainly of theoretical value.

\bigskip

To simplify the problem somewhat, consider the question of whether $C(n)>0$ for a fixed but arbitrary $n$. This is a decision problem, admitting only ``yes'' or ``no'' as an answer. Today, the only certain answers that can be given to this question are positive answers, which can be given only at those $n$ where Costas arrays have been constructed or found through exhaustive enumeration. In particular, it cannot be positively said for any $n>1$ that $C(n)=0$, although there exist plenty of $n$ for which no Costas array is currently known, the smallest two being $n=32$ and 33 (see, for example, \cite{D10,DIR,GT}, and virtually any publication on Costas arrays). Arguing at a high level and in computer-scientific terms, computing the answer to this decision problem involves some computational complexity expressed in terms of $n$. At present, this is extremely high, essentially $O(n!)$, as the only procedure known to yield a guaranteed answer is exhaustive search, whereby every permutation of order $n$ is tested for the Costas property \cite{DIR,DRBCIOW}; though there are, of course, shortcuts, the computation remains of exponential complexity. Can one hope for better? Assuming a candidate for a Costas permutation of order $n$ is provided by an oracle, the actual verification of the Costas property is very fast, involving $O(n^3)$ comparisons \cite{BDR,SVM}, namely of polynomial complexity. These two facts prove that this decision problem lies in NP \cite{GJ}. Is this the best that can be hoped for?

\begin{pr}
Let $E(n)$ stand for the decision problem ``Is $C(n)>0$?''. What is the computational complexity of $E(n)$ expressed in terms of $n$? In particular, is $E$ NP-complete?
\end{pr}

\subsection{Bound on the number of Costas arrays}

Existence of Costas arrays in all orders is equivalent to establishing a positive lower bound for $C(n)$. Of considerable importance would also be an upper bound of $C(n)$. Trivially, $C(n)<n!$, and, unfortunately, the best possible result published today \cite{D10} does not represent a considerable improvement, as it states that
\be \frac{C(n)}{n!}=O\left(\frac{1}{n}\right)\Leftrightarrow C(n)=O((n-1)!).\ee
(It should be noted that the proof of this result, due to D.\ Huw Davies, was presented in the present author's review paper on Costas arrays \cite{D10} and given a nonexisting citation in the literature due to a misunderstanding. It later transpired that D.\ Huw Davies's paper, which the present author was privately handed and consulted for his review \cite{D10}, was, in fact, never published. The proof itself is an excellent instance of the probabilistic method \cite{AS}.)

Though this result is enough to establish that the set of Costas arrays of order $n$ represents asymptotically a set of zero density within the set of all permutations of order $n$, namely that $\lim C(n)/n!\to 0$, the rate of convergence to 0 it gives is too slow compared to what is observed in practice, where $C(n)/n!$, $n\in[29]$ seems to decay to 0 monotonically and at an exponential rate. Two notable efforts towards obtaining an improved bound should be mentioned, both relying on the concept of ``degrees of freedom''.

The first method \cite{SVM} focuses on the entries of the difference triangle of a Costas permutation. After estimating the number of independent comparisons of pairs of entries in the difference triangle that need to be carried out in order to ascertain that a permutation has the Costas property, it estimates the probability of a repeated entry in some row of the difference triangle, and finally estimates the probability $C(n)/n!$ under the simplifying assumption that pairs of repeated entries occur independently. Simplifying slightly, the following formula was reached:
\be \frac{C(n)}{n!}\approx \left(1-\frac{K}{n}\right)^\frac{n^3}{12}\Leftrightarrow C(n)\approx \sqrt{2\pi}\exp\left(-K\frac{n^2}{12}-n+(n+0.5)\ln(n)\right),\ee
using Stirling's approximation. $K$ was selected in \cite{SVM} after fitting the curve to the (then) known values of $C(n)$, $n\in[17]$, and was found to be $K\approx 1$. Amazingly, this estimate for $C(n)$ can easily be verified to remain valid throughout the entire main ``lobe'' of $n\in[26]$ (see Section \ref{kcasec}). On the downside, though, not only is it not rigorous, but it also erroneously predicts that $\lim C(n)=0$.

The second, due to the present author, is inspired by the first, and especially by the underlying principle of ``degrees of freedom'' it introduces, but it looks for them in a different context. More precisely, this concept is now defined \cite{D0} as the smallest $k\leq n$ such that there exist $i_1<i_2<\ldots<i_k\in [n]$ with the property that, for any arbitrary collection of values $f_1,\ldots,f_k\in[n]$, the set of functions $\{f:[n]\to[n]:\ f(i_j)=f_j,\ j\in[k]\}$ contains at most one Costas permutation. Clearly, it follows that
\be C(n)\leq (n-k+1)!\Leftrightarrow \frac{C(n)}{n!}\leq \frac{1}{(n-k)!},\ee
which, depending on the asymptotic behavior of $k=k(n)$, may prove the asymptotic decay of $C(n)/n!$ to 0 at an exponential rate. Indeed, the present author conjectured in \cite{D0} that $k=3$ independent of $n$ for all ``large'' $n$ (more precisely, $n\geq 24$).

In any case, three interesting problems have emerged:

\begin{pr}
Determine the ``degrees of freedom'' of Costas arrays of order $n$, as defined in the second approach above. If possible, improve on this definition.
\end{pr}

\begin{pr}
Determine a bound for $C(n)$ that reflects more closely the actual rate of decay of $C(n)/n!$ observed in practice. Furthermore, prove or disprove that $C(n)/n!$ converges to 0 monotonically.
\end{pr}

\subsection{Sporadic Costas arrays}

Virtually all published results on Costas arrays concern algebraically constructed Costas arrays. Almost nothing is known about sporadic Costas arrays, a fact reflecting the inability to deduce properties of Costas arrays based on the definition alone. The enumeration of orders 28 and 29 revealed no sporadic Costas arrays there, and these are the first orders larger than order 5 where this occurs \cite{DIR,DIRW}. There is, then, the possibility that sporadic Costas arrays cease to exist from a certain order onwards. If so, the situation would be reminiscent of the classification of finite simple groups, which were found to consist of finitely many infinite families plus a finite number of some sporadic groups \cite{W}. There is also the possibility that they constitute examples of generally applicable but still unknown construction techniques.

\begin{pr}
Settle the status of sporadic Costas arrays: are they finitely many? Are there new generally applicable construction techniques that produce them?
\end{pr}

\bigskip

A curious phenomenon associated with sporadic Costas arrays is the existence of twin Costas arrays:

\begin{dfn}
Let $f:[n]\to [n]$ be a Costas permutation, and let $g_1,g_2:[n]\to [n]$ be defined through the relations $g_j(i)=f(i-1)+1,\ i=2,\ldots,n+1,\ j=1,2$, $g_1(1)=g_2(n+2)=1$, and $g_2(1)=g_1(n+2)=n+2$. If $g_1,g_2$ are also Costas permutations, they will be called \emph{twin Costas permutations.}
\end{dfn}

The only known example of twin Costas arrays, as is directly verifiable over the database, is generated out of a sporadic Costas array of order $n=21$.

\begin{pr}
Settle the status of twin Costas arrays: is the only known example the only one in existence? Otherwise, are there finitely/infinitely many pairs of twin Costas arrays? Can they be systematically constructed?
\end{pr}

\bigskip

Another result related to sporadic Costas arrays is yet another question first asked by S.\ Golomb and H.\ Taylor in 1984 \cite{GT}: are there Costas arrays representing configuration of non-attacking queens on the (generalized $n\times n$) chessboard? This question still remains open for $n>1$. However, it is now known \cite{DGR2} that, if such a Costas permutation exists, then it has to be sporadic. It is also known \cite{D2} that such a configuration is constructible on an infinite array. Refining the original question, then,

\begin{pr}\label{queensprob}
Do (sporadic) Costas arrays of finite order $n>1$ exist, which also represent configurations of non-attacking queens on an $n\times n$ chessboard? If yes, how many are they, and can they be systematically constructed?
\end{pr}

There can only be finitely many such arrays: this follows from a surprising connection with honeycomb arrays \cite{BEMP} (more about this in Section \ref{honarsec}).

\bigskip

The final question about sporadic Costas arrays is one for which data is very limited. It has been shown \cite{DGR3} that $G_2$ permutations contain $G_2$ sub-permutations, in the following sense:

\begin{dfn}
Let $f:[n]\to [n]$ be a permutation. If $a_1,a_2\in \NN$ and $b_1,b_2,l\in\NN^*$ exist such that the sets $S_j=\{a_j+b_j i:\ i\in[l]\}$, $j=1,2$ are both subsets of $[n]$ and satisfy $f(S_1)=S_2$, then $g:[l]\to [l]$, $g(i)=(f(a_1+b_1i)-a_2)/b_2$ is said to be a sub-permutation of $f$.
\end{dfn}

\begin{pr}
Is it true that for any Costas permutation $g$ there exists a Costas permutation $f$ such that $g$ is a sub-permutation of $f$? Is this true, in particular, for specific families of Costas permutations (e.g.\ sporadic)?
\end{pr}

\subsection{Cycle structure of Costas permutations}

The iterative composition of $G_2$ permutations has some interesting properties: for example \cite{D0}, in any extension field whose size is a power of a prime number whose exponent contains an odd factor, Lempel $G_2$ permutations are obtainable through the iterative composition of other $G_2$ permutations! This suggests a closer study of the cycles of Costas permutations. Indeed, cycles have normally been part of any classical study of families of permutations, but have not been studied thoroughly in the context of Costas permutations.

\begin{dfn}
Let $f:[n]\to [n]$ be a permutation, and let $I:[n]\to [n]$ be the identity: $I(i)=i$. Let $f^{(0)}=I$, and, for any $k\in\NN$, let $f^{(k+1)}=f(f^{(k)})$, in the sense that, for any $i\in[n]$, $f^{(k+1)}(i)=f(f^{(k)}(i))$.
The cycle of $i$ under $f$ is the set $\{f^{(k)}(i):\ k\in\NN\}$, which is necessarily finite, containing at most $n$ elements.
\end{dfn}

As an example, let us determine the cycles of the $W_1$ permutation $f=[1,2,4,8,5,10,9,7,3,6]$: clearly $f(i)=i$ for $i=1,2,5$, so each of these three points forms its own cycle with a single element. Further, $f(3)=4, f(4)=8, f(8)=7, f(7)=9, f(9)=3$, and $f(6)=10, f(10)=6$, leading to the cycles $(3,4,8,7,9)$ and $(6,10)$. Consequently, we can write $f=(1)(2)(5)(3,4,8,7,9)(6,10)$ to show the cycle structure of $f$ explicitly. Note that the order in which the cycles appear is immaterial; furthermore, the elements within a cycle can be cyclically shifted at will.

Consider now $f^{(2)}$. Given the cycle notation of $f$, it is easy to determine it values: for any $i$, we simply locate its cycle, and the value $f^{(2)}(i)$ will be the element of the cycle lying two positions to the right of $i$, for example $f^{(2)}(7)=3$. More generally, for any $k$, the value $f^{(k)}(i)$ will be the element of the cycle lying $k$ positions to the right of $i$. It follows that $f^{(k)}=I$ iff $k$ is a multiple of the least common multiple of the cycle lengths of $f$, which, in this case, is $2\cdot 5=10$.

\begin{pr}
Describe the cycle structure of Costas permutations. In particular,
\begin{itemize}
  \item Determine the cycle structure of a particular Costas permutation;
  \item Determine the existence/absence of cycles of a particular length in a specific Costas permutation, or a family thereof;
  \item Determine the longest cycle in a specific Costas permutation, or a family thereof;
  \item Determine the smallest $k$ for which $f^{(k)}=I$ for a specific Costas permutation $f$, or a family thereof.
\end{itemize}
\end{pr}

\section{Problems related to the algebraic construction techniques}

\subsection{Disjointness of the algebraic construction techniques}

The large number of algebraic constructions available for Costas arrays begs the question whether they all are genuinely distinct: can one rely on the fact that the sets of Costas arrays produced by each method do not overlap, at least for large enough orders? Or is it possible, for example, for a Costas array to be simultaneously $W_i$ and $G_i$, for some $i\in [4]-1$, or even, assuming $p$ and $p+2$ are twin primes, to be simultaneously $W_1(p)$ and $G_3(p+2)$? Such questions were addressed in \cite{DGR1}, where a comparison between all possible pairs of methods was attempted: most pairs were shown to be incompatible, at least for large enough orders, but for some pairs it was not possible to reach a conclusion.

Analysis of the Costas arrays database suggests that all methods discussed in Section \ref{conalgsec} produce distinct Costas arrays for orders $n>5$. Hence, the question remains:

\begin{pr}
Are the construction methods discussed in Section \ref{conalgsec} pairwise disjoint for orders $n>5$?
\end{pr}

\subsection{Scope of ``dot addition'' methods}

All unpredictably emergent construction techniques described in Section \ref{conalgsec} rely on the addition of a new dot (usually, but not always, a corner dot) to and already existing Costas array. In ``low'' orders, these methods have scored some major successes: for example, the only two known Costas arrays ECs of order 19 prior to their enumeration, and the only one of order 31, are $W_0$ arrays. Moreover, the Rickard construction techniques yielded two new $RW_0$ ECs of order 29, and one new $RG_1$ EC in orders 36 and 42, respectively \cite{R2}.

In higher orders, however, the story is quite different: the Rickard constructions fail to yield any further new Costas arrays up to order 100 \cite{R2}, and the present author verified that this is the case up to order 300. Moreover, in the same search $G_1$ and $W_0$ arrays were also considered, and none were found above orders 52 and 53, respectively (the discovery of the latter was published in \cite{YOH}).

\begin{pr}
Determine the scope of the unpredictably emergent construction techniques discussed in Section \ref{conalgsec}. In particular, is it the case that, for each method amongst them, a $N\in\NN$ can be found, such that this method produces no Costas arrays of order $n>N$?
\end{pr}

\subsection{Alternative construction techniques}\label{actsec}

Considerable effort has been expended towards the discovery of construction techniques for Costas arrays. This area is undoubtedly dominated by the algebraic construction techniques discussed in Section \ref{conalgsec}, and the difficulty associated with discovering new constructions can perhaps be appreciated by the fact that several other proposed techniques either fail to work or are equivalent to the existing ones. For example, J.K.\ Beard's spin-add technique \cite{B2} or polynomial generators \cite{B} are equivalent to the unpredictably emergent methods, while Popovic's construction \cite{P} is, in fact, a rediscovery of the Logarithmic Welch Costas arrays.

An unfortunate side-effect of the long and honored academic tradition of publishing successful research is that researchers are doomed to try out unsuccessful ideas time and again. It is safe to say that it is not possible to know the full range of candidates for Costas arrays generation methods that have been tested and failed, though a published record of them would be invaluable, as it would permit researchers not to repeat mistakes of the past.

\bigskip

For example, let us consider the simple idea of ``interlacing'' two Costas arrays:

\begin{dfn}
Let $A^1=[a^1_{ij}]$ be an $n\times n$ array, and let $A^2=[a^2_{ij}]$ be either an $n\times n$ or an be an $(n-1)\times (n-1)$ array. Let $A=[a_{ij}]$ be the array defined by
\[ a_{ij}=
\begin{cases}
a^1_{\frac{i+1}{2},\frac{j+1}{2}},&  i\bmod2=j\bmod2=1;\\
a^2_{\frac{i}{2},\frac{j}{2}},&  i\bmod2=j\bmod2=0;\\
0, &\text{otherwise}.
\end{cases}
\]
$A$ is said to result from interlacing $A_1$ and $A_2$.
\end{dfn}

Unfortunately, interlacing Costas arrays fails to produce a new Costas array as long as the order of either array is 3 or above: the case of equisized $A_1$ and $A_2$ was treated already in 1985 \cite{FL}, while the case of sizes differing by 1 was treated in \cite{DRG}. The idea is that when size is large enough, both arrays contain a common distance vector. As interlacing rescales distance vectors by a factor of 2 without distorting the arrays, a parallelogram is formed.

\bigskip

An even simpler obvious candidate for a construction technique is concatenation: given $A_1$ and $A_2$ Costas arrays, is
\[A=\left[\begin{array}{cc} A_1 & 0\\ 0 & A_2\end{array}\right]\]
a Costas array? Unfortunately, the argument used in \cite{DRG,FL} is no longer applicable here, since $A_1$ and $A_2$ may be of arbitrary sizes. A stronger result would be needed, namely that any two ``large enough'' Costas arrays have a common distance vector. This seems likely to be the case, but no rigorous proof exists today.

Note that, assuming (without loss of generality) that $A_1$ is no larger than $A_2$, concatenation works for infinitely many orders if $A_1$ is $1\times 1$ or $2\times 2$ (examples are $W_1$ arrays with $c=0$ and $G_2$ arrays with $\alpha+\beta=1$, especially in fields of characteristic 2). However, the largest concatenated array in the database for which $A_1$ is $3\times 3$ is a $7\times 7$ Costas array.

\begin{pr}
Is it the case that any two ``large enough'' Costas arrays have a distance vector in common? It may suffice that the smallest of the two be at least $3\times 3$ and the largest at least $5\times 5$.
\end{pr}

Note that this question was studied in \cite{DGR2}, but only for algebraically constructed Costas arrays.

\bigskip

Having considered a method that fails to deliver, and a method that most likely fails to deliver, let us end with a completely different approach that still looks promising, namely a stochastic search for Costas arrays using genetic algorithms. The underlying idea is to start with a random permutation and, through a sequence of (possibly random) mutations, eventually reach a Costas permutation.

A simple way to implement this search is through a ``gradient descent'' approach: after defining a metric/weight to measure how much a permutation deviates from the Costas property, a mutation is applied on it (from within a pre-specified set of permissible mutations) so as to decrease this metric (ideally, as much as possible). If this can happen under many different mutations, one is applied at random. If no mutation can achieve this, a random mutation is applied, or the procedure restarts with a new random permutation. The procedure runs until a Costas permutation is reached.

A simple such metric is the total sum of the number of repeated entries within each row of the difference triangle of the permutation: for example, if a row contains three equal entries, this contributes two repetitions to the sum. Permissible mutations can be all possible rearrangements of the values of the permutation so that at least a certain pre-specified number of values do not get reassigned: for example, assuming the order is $n$, demanding that at least $n-2$  values do not get reassigned is equivalent to considering all pairwise swaps; demanding that at least $n-3$  values do not get reassigned is equivalent to considering all 3-tuple swaps (which includes pairwise swaps as a subset); and so on. Alternatively, each recorded repetition in the difference triangle can be associated with the three or four values of the permutation giving rise to it, and then the values showing the most such ``hits'' can be deranged amongst themselves.

Stochastic searches were considered in \cite{DHR,Dc1,R7,RH}, in which several different algorithms were tested. The general conclusion was that, although they performed adequately for ``low'' orders, they consistently failed to recover any Costas arrays in ``higher'' orders within a reasonable time interval (more specifically, in orders above 15, allowing for a runtime up to 1 hour). As is the case with all genetic algorithms, success is heavily dependent upon the selection of a suitable metric and set of permissible mutations.

\begin{pr}
Develop a stochastic search method for Costas arrays which works in large enough orders of interest (in particular, in orders 30, 31, 32, and 33).
\end{pr}

\section{Inverse problems}\label{ipsec}

Most publications discussing properties of Costas arrays prove statements of the form: ``If a Costas array is in family $F$, then it has property $P$''. For example, a $T_4$ array represents a configuration of non-attacking kings on a chessboard \cite{DGR2}; a $W_1$ array is singly periodic; a $W_1$ array has anti-reflective symmetry \cite{BCGRMP}; the special subfamily of symmetric $G_2$ permutations of order $r^2-2$ discussed in \cite{DGO} have $r$ fixed points; and so on. There are very few statements, if any, of the form: ``If a Costas array has property $P$, then it must be a member of family $F$.'' These problems will be referred to as inverse problems, and they have proved extremely difficult to treat so far.

\subsection{Single periodicity}

A quick search in the database identifies sporadic arrays possessing anti-reflective symmetry, so this property does not characterize $W_1$ arrays. However, no singly periodic Costas array other than $W_1$ arrays is known; it is only known that a Costas array of odd order cannot be singly periodic \cite{GT}. Quoting then a problem out of \cite{GT} once more:

\begin{pr}
Is it true that every singly periodic Costas array is a $W_1$ array?
\end{pr}

A significant step towards the solution of this problem was taken in \cite{GM}, where the concept of a circular Costas array was introduced:
\begin{dfn}
A permutation $f:[n]-1\to [n]-1$ is said to be \emph{circular Costas}, for all $k\in[n-1]$, all vectors $f(i+k)-f(i)$, $i\in[n]-1$ are distinct $\bmod (n+1)$, when $i+k$ is considered $\bmod n$. 
\end{dfn}

Clearly, every $W_1$ permutation is circular Costas, but the converse may not be true. It was shown in \cite{GM} that, if a circular Costas permutation of order $n$ exists, then $n+1$ must be a prime. 

\subsection{Cyclic shifts when blank columns/rows are added}

It was mentioned earlier, in Section \ref{conalgsec} and in connection with the Rickard construction \cite{R2}, that $W_1^\text{exp}$ arrays have the property that, if a blank row is added at the top and then rows are cyclically shifted, the resulting array, which is no longer square, still has the Costas property; and similarly, that $G_2$ arrays have the property that, if a blank row and a blank column are added at the top and at the side, and then rows and columns are cyclically shifted, the resulting array, which is square but no longer represents a permutation, still has the Costas property. These properties are direct consequences of the definitions. Are, however, the converse statements true?

\begin{pr}
Let a Costas array have the property that, if a blank row is added at the top, then any cyclic shift of the rows leads to a (non-square) array with the Costas property. Is such a Costas array necessarily a $W_1^\text{exp}$ array?
\end{pr}

\begin{pr}
Let a Costas array have the property that, if a blank row and a blank column is added at the top and at the side, respectively, then any cyclic shift of the rows and columns leads to a (square but non-permutation) array with the Costas property. Is such a Costas array necessarily a $G_2$ array?
\end{pr}

%

\section{Problems related to properties of Costas arrays}

\subsection{Parity populations of $G_2$ arrays in characteristic 2}

It was mentioned in Section \ref{actsec} that interlacing two Costas arrays cannot lead to a new Costas array, as long as the order of either array is 3 or above. An interlaced array is considerably constrained, as its dots can only lie in positions whose coordinates are both even or both odd (this is a direct consequence of the definition). How close do algebraically constructed Costas arrays get to this state? Let $ee$, $eo$, $oe$, and $oo$ stand for the number of dots of the Costas array whose coordinates are both even, even and odd, odd and even, and both odd, respectively: it was found in \cite{DGR3} that, for $G_2$ and $W_1$ arrays, these four quantities, called the \emph{parity populations}, were not only as equal as they could be (with the significant exception of $W_1$ arrays constructed in $\FF(p)$, $p\equiv 3\bmod 4$, where they were found to depend on the class number), but they were also independent of the primitive root(s) used. The latter was no longer the case, however, for $G_2$ arrays constructed in fields of characteristic 2, where the parity populations varied across the various $G_2$ arrays in the family.

\begin{pr}
Determine the parity populations of $G_2$ arrays constructed in fields of characteristic 2.
\end{pr}

Some numerical results relevant to this problem can be found in \cite{D9}.

\subsection{Fixed points}

It was mentioned in Section \ref{ipsec} that a special subfamily of symmetric $G_2$ permutations of order $q-2$, $q$ a square prime power, have asymptotically $\sqrt{q}$ dots on their main diagonal. This is an interesting result, as diagonals of Costas arrays are Golomb rulers \cite{D7}, and it is conjectured that the largest number of dots a Golomb ruler of a given length can contain is asymptotically equal to the square root of its length \cite{D7}: consequently, this subfamily is asymptotically optimally dense. Unfortunately, this turns out not to be a new construction technique for Golomb rulers, but rather equivalent to a classical one \cite{D7}.

This result motivates the study of the number of fixed points of other Costas permutations. For example, setting $f(i)=i$ and $\beta=\alpha^r$ for some $r$ such that $(r,q-1)=1$ in Theorem \ref{g2}, and then setting $\alpha^i=x$, it follows that the number of fixed points of this $G_2$ permutation equals the number of roots of $x^r+x=1$ in $\FF(q)$, excepting $x=1$. Far more interesting, though, is the case of $W_1$ permutations: setting $f(i)=i$ in Theorem \ref{w1} it follows that the number of fixed points of a $W_1$ permutation equals the number of roots $i$ of
\be i\equiv \alpha^{i-1+c}\bmod p\Leftrightarrow i\equiv C\alpha^i\bmod p,\ C=\alpha^{c-1}.\ee

This is a transcendental equation over a finite field, and, to the best of our knowledge, there is no precedent of any relevant published work in the literature. The closest question considered is whether, for a given $i$, there exists a primitive root $\alpha$ such that $i\equiv \alpha^i\bmod p$, which is, in effect, the inverse problem, and has been answered in the affirmative (see \cite{HM,LPS,Z} and references therein).

\begin{pr}
Let $q$ be a prime power and $r$ be such that $(r,q-1)=1$. Determine the number of roots $x$ of the equation $x^r+x=1$ in $\FF(q)$.
\end{pr}

\begin{pr}
Let $p$ be a prime, $C\in \FF^*(p)$, and $\alpha$ a primitive root of $\FF(p)$. Determine the number of roots $i \in [p-1]$ of $i\equiv C\alpha^i\bmod p$.
\end{pr}

Some numerical results relevant to this problem can be found in \cite{D9}. Moreover, a stochastic model relating the solution of this problem to Lambert's function can be found in \cite{D8}.

\subsection{Forbidden positions}

Consider overlaying all known Costas arrays of order $n$. Are there any positions of the $n\times n$ square grid not covered by a dot? These positions, called forbidden positions, can simplify the enumeration of Costas arrays whenever they exist, as any permutation array featuring a dot in a forbidden position is automatically disqualified from being Costas. A search over the database proves that there is a single forbidden position in order 3 (namely the midpoint (2,2) of the grid), that there are no forbidden positions in orders $4\leq n\leq 24$, and that forbidden positions appear again in orders 25, 26, 27, and 29 \cite{DIR,DIRW} (clearly, there can be no forbidden positions in orders of the form $p-1$, $p$ prime, due to the single periodicity of the $W_1$ construction which applies there).

Unfortunately, forbidden positions have so far only been determined after the complete enumeration of a certain order, not before, with the exception of order 3. Though this order is trivially small, there is a brief and elegant argument that can be used, which may be generalizable to higher orders. Indeed, the first row of the difference triangle of a permutation $f$ of order 3 contains only two entries, namely $f(2)-f(1)$ and $f(3)-f(2)$, which are equal iff
\[f(2)-f(1)=f(3)-f(2)\Leftrightarrow 2f(2)=f(1)+f(3)\Leftrightarrow 3f(2)=f(1)+f(2)+f(3)=6\Leftrightarrow f(2)=2.\]
The key fact in this argument is that symmetric (polynomial) expressions over the values of $f$ are the same for all permutations $f$: in this case, $f(1)+f(2)+f(3)=1+2+3=6$. Can this technique be used for higher orders as well?

\begin{pr}
Determine the forbidden positions for Costas arrays of a certain order without resorting to enumeration.
\end{pr}

\subsection{Cross-correlation of $G_2$ and $W_1$ arrays}

Cross-correlation is one of the aspects of Costas arrays that has attracted considerable attention, in view of its importance in applications. More specifically, the optimal suppression of sidelobes exhibited by waveforms described by Costas arrays is, in practice, only half of the waveform's desirable behavior. As many RADAR/SONAR systems may be operating within the same geographical area and in the same frequency spectrum, provisions should be made to minimize cross-interference resulting from waveform collisions \cite{DT,MST}. The obvious solution is to select Costas arrays from within a family where pairwise cross-correlation of Costas arrays is uniformly bounded by as low as possible a bound. It would not be an exaggeration to say, then, that low cross-correlation is the ``second half'' of the definition of Costas arrays that find their way in practical RADAR/SONAR systems applications.

Perhaps the first result relevant to the cross-correlation of Costas arrays was the demonstration by A.\ Freedman and N.\ Levanon \cite{FL} in 1985 that any two Costas arrays of the same order (trivially small orders excepted) must have a maximal cross-correlation of at least 2, and hence that ideal cross-correlation and ideal auto-correlation are incompatible properties. The study of cross-correlations became more focused in later years, with the work by W.\ Chang and K.\ Scarbrough in 1989 \cite{CS}, which demonstrated that there exist sets of Welch Costas arrays whose pairwise cross-correlation is ideal, provided that only a strip of the correlation surface around the origin is considered, and not the entire surface. Later, in 1994, S.\ Maric, I.\ Seskar, and E.\ Titlebaum \cite{MST} studied the cross-correlation of Welch Costas arrays more thoroughly, providing a table with the maximal observed cross-correlation between any two Welch Costas arrays of order $p-1$, $p$ prime (effectively a first version of our Table \ref{dat} below), and showing that this cross-correlation is bounded above by $\ds \frac{p-1}{2}$. The conference version of this work had appeared in ICASSP 1990 \cite{TM}, where the authors had already observed that some primes corresponded to local minima in the table, and had identified them to be safe primes.

An important part of the puzzle was still missing, though, namely cross-correlation results for Golomb Costas arrays. The subject had been considered by D.\ Drumheller and E.\ Titlebaum in 1991 \cite{DT}, who suggested an upper bound that is, unfortunately, not always tight. Then, S.\ Rickard, in his Master's thesis in 1993 \cite{R}, extended the cross-correlation tables for Welch arrays of \cite{TM} and offered the corresponding results for Golomb arrays. He observed a surprising link between the Welch and Golomb results (given the two construction methods are so dissimilar), confirmed the special role of the safe primes for the Golomb results as well, and asked whether this link could be explained, though he did not offer any explanation.

A unified presentation and formal justification of all aspects of cross-correlation of Costas arrays mentioned above was attempted in \cite{DGRST}. Though this publication went indeed considerably beyond its predecessors, it too did not succeed in explaining all observed results: many of them are listed in \cite{DGRST} as conjectures supported by empirical evidence and awaiting formal proof.

\begin{pr}
Let $p$ stand for a prime and $q$ for a prime power. Define $\ds \Psi_{W_1}(p)=\max_{(u,v)}\max_{\left(f_1\overset{\text{pr}}\neq f_2\right)}\Psi_{f_1,f_2}(u,v)$ and $\ds \Psi_{G_2}(q)=\max_{(u,v)}\max_{(f'_1\neq f'_2)}\Psi_{f'_1,f'_2}(u,v)$, where the notation $\ds f_1\overset{\text{pr}}\neq f_2$ denotes that $f_1$ and $f_2$ must be produced by different primitive roots.

\medskip

Prove that, with respect to $\Psi_{W_1}(p)$ and $\Psi_{G_2}(p)$, primes $p$ can be classified into three groups:
\begin{itemize}
  \item For non-safe primes $p\neq 19$ (such that $\ds \frac{p-1}{2}$ is not a prime, that is), $\ds \Psi_{W_1}(p)=\Psi_{G_2}(p)+1=\frac{p-1}{t}$, where $t$ is the smallest prime such that $p\equiv 1\bmod(2t)$.
  \item $p=19$ is the only prime for which $\ds \Psi_{W_1}(p)=\Psi_{G_2}(p)=\frac{p-1}{t}=6=\frac{19-1}{3}$, 3 being the smallest prime $t$ such that $19\equiv 1\bmod(2t)$.
  \item Primes $p$ such that $\ds \frac{p-1}{2}$ is also a prime (safe primes, that is) correspond to local minima of both $\Psi_{W_1}(p)$ and $\Psi_{G_2}(p)$.
\end{itemize}
Produce a formula for the values of these local minima.

\medskip

Furthermore, prove that, with respect to $\Psi_{G_2}(q)$, prime powers $q$ can be classified into three groups:
\begin{itemize}
  \item For prime powers $q\neq 16$ such that neither $\ds \frac{q-1}{2}$ nor $q-1$ is a prime, $\ds \Psi_{W_1}(q)=\frac{q-1}{t}-1$, where $t$ is either the smallest prime such that $q\equiv 1\bmod (2t)$, if $q$ is odd, or such that $q\equiv 1\bmod t$, if $q$ is even.
  \item $q=16$ is the only prime power for which $\ds \Psi_{G_2}(q)=\frac{q-1}{t}=5$, where 3 is the smallest prime $t$ such that $16\equiv 1 \bmod t$.
  \item Prime powers $q$ such that either $\ds \frac{q-1}{2}$ is prime (in which case, as can be easily verified, $q=3^m$ for some prime $m$) or $q-1$ is prime (in which case $q=2^m$ for some prime $m$ and $q-1$ is a Mersenne prime) correspond to local minima of $\Psi_{G_2}(q)$.
\end{itemize}
Produce a formula for the values of these local minima.

\medskip

Assuming $p$ is a non-safe prime, show that $\ds \Psi_{W_1}(p)=\max_{(f_1\overset{\text{pr}}\neq f_2)}\Psi_{f_1,f_2}(0,0)$. Assuming $q\neq 16,19$ is a non-safe prime power such that $q-1$ is not prime, show that $\ds \Psi_{G_2}(q)=\Psi_{f'_1,f'_2}(0,0)$, where $f'_1$ and $f'_2$ are $G_2$ permutations generated in $\FF(q)$, such that, assuming the former is generated by primitive roots $\alpha$ and $\beta$, the latter is generated by $\alpha^r$ and $\beta$, where $r$ is defined through the requirements that a) $\ds r=\lambda\frac{p-1}{w}+1$; b) $(r,q-1)=1$; c) $\lambda=1$ or 2; and d) $w$ is the smallest possible (such an $r$ is guaranteed to exit).

\end{pr}

Full details regarding the formulation of this detailed and complex conjecture can be found in \cite{DGRST}. The cross-correlation of Costas arrays within the same EC was considered in \cite{DGHR}.

The consideration of cross-correlation suggests that not all Costas arrays may be equally suitable for RADAR/SONAR applications. Another criterion to consider is the maximal hop between consecutive dots/frequencies in a Costas array, as the farther apart two frequencies lie, the more the response of frequency hopping filter deviates from the ideal response when hopping between them \cite{D5}.

\subsection{A necessary property for $G_2$-arrays}

It was mentioned above that a necessary property of $W_1$ arrays is their anti-reflective symmetry: every $W_1$ array is anti-reflective, but there are some sporadic Costas arrays that are anti-reflective as well. A necessary property for $G_2$ arrays generated in finite fields of odd characteristic was published in \cite{D1}: letting $q$ be the size of the field, it follows that, for any $G_2$ permutation $f$,
\be\label{g2necpr} f(\mu+i)-f(\mu-i)\equiv i[f(\mu+1)-f(\mu-1)]\bmod (q-1),\text{ where } \mu=\frac{q-1}{2}.\ee
Verification over the database suggest that this property may also be sufficient for $G_2$ arrays, as no non-$G_2$ Costas array was found to satisfy it.

\begin{pr}
Find a necessary property, similar in nature to anti-reflective symmetry and (\ref{g2necpr}), for $G_2$ arrays generated in finite fields of characteristic 2.
\end{pr}

\section{Problems related to variants and generalizations of Costas arrays}

Costas arrays are predominantly characterized by their optimal autocorrelation property, and as such they are only one out of the many different mathematical objects defined around this particular principle (each in its own distinct context) and considered in the literature. Examples include Sidon sets (also known as Golomb rulers) \cite{D7,R4}, RADAR/SONAR sequences (also known as Golomb rectangles) \cite{BT,G0,HZ,MGC,MGT,R4,R5,ZT}, Golomb rectangles \cite{R5}, honeycomb arrays \cite{BEMP,BPPS,GT}, as well as arrays of dots where ``half'' of the Costas property hold, namely where all linear segments connecting pairs of dots either have distinct lengths or distinct slopes (but not necessarily both) \cite{EGRT,LT,PT,Z2}. Further examples include generalizations of the Costas property in higher dimensions \cite{D4,E} and in the continuum \cite{D3,D6,DR2}. Some Costas arrays were also found found to yield Almost Perfect Nonlinear (APN) permutations \cite{DGM,DMR}, which find applications in cryptography. An extensive study of signals built for ``good correlation'' (in various contexts) is presented in \cite{GG2}.

\subsection{Golomb rectangles}

Progressive relaxation of the constraints present in the definition of Costas arrays leads to new mathematical objects. For example, binary $m\times n$ arrays such that there be exactly one dot per column, and such that the range of their (horizontal only) autocorrelation consists only of the three values $\{n,1,0\}$ are known as RADAR/SONAR arrays/sequences \cite{G0,MGT,R4} and find applications in RADAR/SONAR systems. One step further, in extended RADAR sequences blank columns are also allowed \cite{MGC}.

The next natural step in this process is to abolish any restriction on the number and position of the dots in the array other than the distinct differences property itself.

\begin{dfn}
A Golomb rectangle \cite{R5} is a rectangular array of $N>2$ dots and blanks such that its auto-correlation range consists of exactly three values: $\{N,1,0\}$.
\end{dfn}

Such configurations also have found important application in cryptography, and specifically in the distribution of cryptographic keys over sensor networks \cite{BEMP2}. Nevertheless, virtually nothing is known about the mathematical properties of these objects.

\begin{pr}
Let $C(m,n,N)$ stand for the number of $m\times n$ Golomb rectangles containing $N$ dots. Calculate $C(m,n,N)$, and, in particular, for given $m$ and $n$, find the maximal $N$ for which $C(m,n,N)>0$. Conversely, for a given $N$, find the smallest values of $m$ and $n$ for which $C(m,n,N)>0$.
\end{pr}

The latter question has received considerable attention in the literature in connection with (optimal) RADAR/SONAR arrays \cite{BT,HZ,MGT,ZT}: for an $m\times n$ RADAR/SONAR array (which has $n$ dots, one per column), what is the maximal $n$ for a given $m$ (or effective bounds thereof)?

\subsection{Honeycomb arrays}\label{honarsec}

A honeycomb array \cite{BPPS,GT} is the equivalent of a Costas array on a hexagonal grid. Ordinary rectangular arrays are implicitly understood as conforming to a square grid: elements are neatly arranged in rows and columns, and each element (except for those at the border of the array) has four natural neighbors, namely those lying above, below, to the left, and to the right, as if a square frame is drawn around each element and elements whose frames share an edge are defined to be neighbors. In a hexagonal grid, on the other hand, the frame around each element is a regular hexagon, so there are three natural directions associated with the grid, and six neighbors to each element (except for those at the border). The side of the regular hexagon can be assumed to be of length $1/\sqrt{3}$, so that centers of neighboring hexagons are at distance 1 apart.

The definition now is straightforward:

\begin{dfn}
Let $r\in\NN$, and consider a hexagonal array of hexagonal cells, whose elements are dots and blanks, such that its side consist of $r$ hexagons (and its diameter by $2r+1$ hexagons), and such that there be exactly one dot per row in each natural direction of the grid (that is, a total of $2r+1$ dots). This will be a honeycomb array iff either of the following equivalent conditions holds:
\begin{itemize}
  \item no two linear segments connecting pairs of dots (assumed to lie at the center of the cell they lie in) have both the same length and the same slope; or
  \item assuming two copies of the array are placed on top of each other, and one is shifted with respect to the other to any position (but not rotated, and so that overlapping cells overlap perfectly and not in part), the number of pairs of overlapping dots can only be 0, 1, or $2r+1$ (this is a restatement of the auto-correlation Definition \ref{ccdefd} in the hexagonal setting).
\end{itemize}
\end{dfn}

Fortunately, there is a bijective correspondence between the hexagonal and the square grid \cite{BEMP,BPPS,GT}, which allows the definition of honeycomb arrays as a subfamily of Costas arrays. Indeed, position a honeycomb array so that its diameter is horizontal, and align rows on top of the diameter along the right (to have a common start with the diameter), and rows below the diameter along the left (to have a common end with the diameter). In this arrangement, two of the three natural directions of the honeycomb array now correspond to rows and columns, while the third corresponds to diagonals along the bottom left--top right direction. Naturally, this is not a complete square array yet, as some square cells on the top left and bottom right corner are missing, but they can be added containing blanks. Thus, a honeycomb array of radius $r$ has been mapped into a Costas array of order $2r+1$, and the mapping can be inverted.

\begin{dfn}
A honeycomb array is a Costas array with the extra property that every diagonal along the bottom left--top right direction contains exactly one dot.
\end{dfn}

A major and surprising result is that no honeycomb array with $2r+1\geq 1289$ dots exists \cite{BEMP} (or, equivalently, of radius $r\geq 644$), and, consequently, only finitely many honeycomb arrays exist. Specifically, only 12 (up to symmetry) are known today, and it is conjectured that no other exists \cite{BPPS}.

\begin{pr}
Enumerate all honeycomb arrays.
\end{pr}

Note that a Costas array configuration of non-attacking queens on the chessboard (see Problem \ref{queensprob}) is a special case of a honeycomb array, where each diagonal along the perpendicular direction, namely bottom right--top left, also contains exactly one dot. Thus, there can be at most finitely many such Costas arrays as well.

\subsection{Relation between Costas arrays and Golomb rulers}

Golomb rulers owe their name to the studies of Prof.\ S.\ Golomb, who demonstrated an important application in graph labeling \cite{G3}. They have important applications in engineering, for example in radio-frequency allocation for avoiding third-order interference \cite{B0}, in generating convolutional self-orthogonal codes \cite{RB}, in the formation of optimal linear telescope arrays in radio-astronomy \cite{BBR} etc., as well as an inherent mathematical interest. Even earlier, they had also been studied in the context of Fourier series \cite{S2}.

A Golomb ruler (or Sidon set, or $B_2$-sequence) is a collection of distinct integers such that all pairwise differences between them are also distinct. It can be visualized as a vector of dots and blanks, so that an element of the vector is a dot iff its position corresponds to one of the selected integers. Clearly, overlaying two copies of the vector, sliding one of the copies, and counting the number of pairs of overlapping dots leads to a count of either 0, 1, or of the total number of dots, if the two copies overlap perfectly. Hence, Golomb rulers are 1D analogs of Costas arrays. Without loss of generality, it can be assumed that the Golomb ruler contains $m$ integers, the smallest of which is $0$ and the largest $n$.

\begin{dfn}
Let $m\in\NN^*$, let $f:[m]\to \NN$ be strictly increasing, and let $f(1)=0$ and $n=f(m)$. The range of $f$ is a \emph{Sidon set} (or a \emph{$B_2$-sequence}) iff all differences $f(i)-f(j),\ 0\leq j<i\leq n$ are distinct. In that case, the vector $v:[n+1]-1\to \{0,1\}$ such that $v(f(i))=1$, $i\in[m]$ and 0 elsewhere is a \emph{Golomb ruler.} $n$ is the \emph{length} of the Golomb ruler, and $m$ its \emph{number of points.}
\end{dfn}

The two fundamental problems in the field are to determine a) the maximal possible $m$ for a fixed $n$, and b) the minimal possible $n$ for a fixed $m$. Golomb rulers satisfying either of these two extremes are called \emph{optimal} (and, more specifically, \emph{optimally dense} in the former case, and \emph{optimally short} in the latter).

\begin{pr}
Provide exact formulas for the length of optimally short Golomb rulers of a given number of points, and for the number of points of optimally dense Golomb rulers of a given length.
\end{pr}

Several inequalities involving $m$ and $n$ have been proposed for optimal Golomb rulers \cite{D,ET,L}: they essentially state that $m\approx \sqrt{n}\Leftrightarrow n\approx m^2$. Furthermore, several construction techniques are available for Golomb rulers, which work for infinitely many (but not all) values of $m$, namely the Erd\"os-Turan \cite{ET}, Ruzsa-Lindstr\"om \cite{L,L2,R6}, Bose-Chowla \cite{B3,BC}, and Singer \cite{S1} constructions. A very comprehensive online resource about Golomb rulers is J.B.\ Shearer's webpage on Golomb rulers \cite{S3}, where all optimal known optimal Golomb rulers are posted, while an overview of the theory of Golomb rulers and their construction techniques is provided in A.\ Dimitromanolakis's Diploma Thesis \cite{D}, and also in \cite{D7}. It should be mentioned that considerable effort has been put into the enumeration of optimal Golomb rulers, just like in Costas arrays (see, for example, \cite{DRM,R8,ogr}). 

Efforts have been made to relate Golomb rulers and Costas arrays directly, by converting one object into the other. A procedure was described in \cite{DR1,R4} whereby any $n\times n$ Costas array (or a Golomb rectangle, more generally) can be converted into a Golomb ruler by adding at least $n-2$ blank rows at the bottom and then stacking its columns in order one below the other. It was also determined empirically in \cite{DR1} that for some Costas arrays less than $n-2$ blank rows could be added with the same effect. A more elegant approach which involves tiling the integer plane using a binary array possessing distinct differences and then applying an unfolding operator is described in \cite{E2}.

\begin{pr}\label{lastprob}
Convert efficiently a Costas array (or Golomb rectangle) into a Golomb ruler and vice versa. In particular, for a given Costas array (or a family thereof), determine the minimal number of blank rows that must be added at the bottom in order to form a Golomb ruler by stacking its columns in order one below the other.
\end{pr}

\section{Conclusion}

After brief reports on recent research progress in various aspects of Costas arrays, a total of 26 open problems were proposed, and specifically:
\begin{itemize}
  \item nine core theoretical problems;
  \item four problems related to the algebraic construction techniques;
  \item three inverse problems;
  \item six problems related to properties; and
  \item four problems on variants and generalizations.
\end{itemize}

This quite diverse selection, which naturally reflects the present author's experience and preferences, is intended to motivate researchers to consider working on the fascinating subject of Costas arrays.

\section*{Acknowledgements}

This material is based upon works supported by the Science Foundation Ireland under Grant No.\ 08/RFP/MTH1164.


\begin{thebibliography}{10}
\bibitem{AS} N.\ Alon and J.\ Spencer. ``The probabilistic method.'' Wiley-Interscience, 2000.
\bibitem{B0} W.C.\ Babcock, Intermodulation interference in radio systems/frequency of occurrence and control by channel selection, Bell System Technical Journal, Volume 31, 1953, pp.\ 63--73.
\bibitem{BDR} L.\ Barker, K.\ Drakakis, and S.\ Rickard. ``On the complexity of the verification of the Costas property.'' Proceedings of the IEEE, Volume 97, Issue 3, March 2009, pp.\ 586--593.
\bibitem{B} J.K.\ Beard. ``Costas array generator polynomials in finite fields.'' Conference on Information Signals and Systems, Princeton University, USA, 2008.
\bibitem{B2} J.K.\ Beard. ``Generating Costas arrays to order 200.'' Conference on Information Signals and Systems, Princeton University, USA, 2006.
\bibitem{BREMW} J.K.\ Beard, J.C.\ Russo, K.G.\ Erickson, M.C.\ Monteleone, and M.T.\ Wright. ``Costas arrays generation and search methodology.'' IEEE Transactions on Aerospace and Electronic Systems, Volume 43, 2007, pp.\ 522--538.
\bibitem{BBR} F.\ Biraud, E.\ Blum, and J.\ Ribes. ``On optimum synthetic linear arrays with application to radioastronomy.'' IEEE Transactions on Antennas and Propagation, Volume 22, Issue 1, 1974, pp.\ 108--109.
\bibitem{BEMP} S.R.\ Blackburn, T.\ Etzion, K.M.\ Martin, and M.B.\ Paterson. ``Two-Dimensional Patterns With Distinct Differences --- Constructions, Bounds, and Maximal Anticodes.'' IEEE Transactions On Information Theory, Volume 56, Issue 3, March 2010, pp.\ 3961--3972.
\bibitem{BEMP2} S.R.\ Blackburn, T.\ Etzion, K.M.\ Martin, and M.B.\ Paterson. ``Distinct Difference Configurations: Multihop Paths and Key Predistribution in Sensor Networks.'' IEEE Transactions On Information Theory, Volume 56, Issue 8, August 2010, pp.\ 1216--1229.
\bibitem{BPPS} S.R. Blackburn, A.\ Panoui, M.B.\ Paterson, and D.R.\ Stinson. ``Honeycomb arrays.'' The Electronic Journal of Combinatorics, Volume 17, 2010.
\bibitem{B3} R.C.\ Bose. ``An affine analogue of Singer's theorem.'' Journal of the Indian Mathematical Society, Volume 6, 1942, pp.\ 1–-15.
\bibitem{BC} R.C.\ Bose and S.\ Chowla. ``Theorems in the additive theory of numbers.'' Commentarii Mathematici Helvetici, Volume 37, 1962-63, pp.\ 141--147.
\bibitem{BCGRMP} C.\ Brown, M.\ Cenki, R.\ Games, J.\ Rushanan, O.\ Moreno, and P.\ Pei. ``New enumeration results for Costas arrays.'' IEEE International Symposium on Information Theory, pp.\ 405, January 1993.
\bibitem{BT} A.\ Blokhuis and H.J.\ Tiersma. ``Bounds for the size of RADAR arrays.'' IEEE Transactions On Information Theory, Volume 34,  Issue 1, pp.\ 164--167, January 1988.
\bibitem{Ch} W.\ Chang. ``A remark on the definition of Costas arrays.'' Proceedings of the IEEE, Volume 75, Issue 4, pp.\ 522--523, April 1987.
\bibitem{CS} W.\ Chang and K.\ Scarbrough. ``Costas arrays with small number of crosscoincidences.'' IEEE Transactions on Aerospace and Electronic Systems, Volume 25, Issue 1, pp.\ 109--113, January 1989.
\bibitem{C} W.\ Correll Jr. ``A closed form expression for the number of Costas arrays of arbitrary order.'' 41st Asilomar Conference on Signals, Systems and Computers, Pacific Grove, CA, USA, November 2008.
\bibitem{C1} J.P.\ Costas. ``Medium constraints on sonar design and performance.'' Technical Report Class 1 Rep. R65EMH33, GE Co., 1965.
\bibitem{C2}  J.P.\ Costas. ``A study of detection waveforms having nearly ideal range-doppler ambiguity properties.'' Proceedings of the IEEE, Volume 72, Issue 8, pp.\ 996--1009, August 1984.
\bibitem{CM} S.\ Cohen and G.\ Mullen. ``Primitive elements in finite fields and Costas arrays.'' Applicable Algebra in Engineering, Communication and Computing, Volume 2, 1991, pp.\ 45--53.
\bibitem{D} A.\ Dimitromanolakis. ``Analysis of the Golomb ruler and the Sidon set problems, and determination of large, near-optimal Golomb rulers.'' Diploma thesis, Department of Electronic and Computer Engineering, Technical University of Crete (2002, available online in English at \texttt{http://www.cs.toronto.edu/~apostol/golomb/}).
\bibitem{DRM} A.\ Dollas, W.T.\ Rankin, and D.\ McCracken. ``A new algorithm for Golomb ruler derivation and proof of the 19 mark ruler.'' IEEE Transactions On Information Theory, Volume 44, Issue 1, January 1998, pp.\ 379--382. 
\bibitem{D0} K.\ Drakakis. ``On the degrees of freedom of Costas permutations and other constraints.'' Advances in Mathematics of Communications.
\bibitem{D1} K.\ Drakakis. ``A structural constraint for Golomb Costas arrays.'' IEEE Transactions on Information Theory, Volume 56, Issue 11, November 2010, pp.\ 5762--5764.
\bibitem{D2} K.\ Drakakis. ``On the existence of infinite size Costas arrays configurations of nonattacking queens on the chessboard.'' International Journal of Combinatorics, Volume 2010.
\bibitem{D3} K.\ Drakakis. ``Cauchy's functional equation and nowhere continuous/everywhere dense Costas bijections in Euclidean spaces.'' (accepted for publication in) ``Functional equations in mathematical analysis --- Editors: J.\ Brzdek and Th.M.\ Rassias.'' Springer.
\bibitem{D4} K.\ Drakakis. ``On the generalization of the Costas property in higher dimensions.'' Advances in Mathematics of Communications, Volume 4, Issue 1, 2010, pp.\ 1--22.
\bibitem{D5} K.\ Drakakis. ``On the hops present in Costas permutations.'' IEEE Transactions on Information Theory, Volume 56, Issue 3, 2010, pp.\ 1271--1277.
\bibitem{D6} K.\ Drakakis. ``On Costas sets and Costas clouds.'' Abstract and Applied Analysis, Volume 2009.
\bibitem{D7} K.\ Drakakis. ``A review of the available construction methods for Golomb rulers.'' Advances in Mathematics of Communications, Volume 3, Issue 3, August 2009, pp.\ 235--250.
\bibitem{D8} K.\ Drakakis. ``A stochastic model for the number of fixed points of a Welch Costas permutation.'' Ars Combinatoria, Volume 92, July 2009, pp.\ 33--52.
\bibitem{D9} K.\ Drakakis. ``Three challenges in Costas arrays.'' Ars Combinatoria, Volume 89, October 2008, pp.\ 167--182.
\bibitem{D10} K.\ Drakakis. ``A review of Costas arrays.'' Journal of Applied Mathematics, Volume 2006.
\bibitem{Dc1} K.\ Drakakis. ``Computer-assisted search for Costas arrays.'' IET Irish Signals and Systems Conference 2005, pp.\ 136--144.
\bibitem{DGHR} K.\ Drakakis, R.\ Gow, J.\ Healy, and S.\ Rickard. ``Cross-correlation properties of Costas arrays and their images under horizontal and vertical flips.'' Mathematical Problems in Engineering, Volume 2008.
\bibitem{DGM} K.\ Drakakis, R.\ Gow, and G.\ McGuire. ``APN Permutations on $\mathbb{Z}_n$ and Costas Arrays.'' Discrete Applied Mathematics, Volume 157, Issue 15, August 2009, pp.\ 3320--3326.
\bibitem{DGO} K.\ Drakakis, R.\ Gow, and L.\ O'Carroll. ``On the symmetry of Welch- and Golomb-constructed Costas arrays.'' Discrete Mathematics, Volume 309, Issue 8, April 2009, pp.\ 2559--2563.
\bibitem{DGR1} K.\ Drakakis, R.\ Gow, and S.\ Rickard. ``On the disjointness of algebraically constructed Costas arrays.'' (accepted for publication in) Journal of Algebra and its Applications.
\bibitem{DGR2} K.\ Drakakis, R.\ Gow, and S.\ Rickard. ``Common distance vectors between Costas arrays.'' Advances in Mathematics of Communications, Volume 3, Issue 1, February 2009, pp.\ 35--52.
\bibitem{DGR3} K.\ Drakakis, R.\ Gow, and S.\ Rickard. ``Parity properties of Costas arrays defined via finite fields.'' Advances in Mathematics of Communications, Volume 1, Issue 3, August 2007, pp.\ 323--332.
\bibitem{DGRST} K.\ Drakakis, R.\ Gow, S.\ Rickard, J.\ Sheekey, and K.\ Taylor. ``On the maximal cross-correlation of algebraically constructed Costas arrays.'' (accepted for publication in) IEEE Transactions on Information Theory.
\bibitem{DHR} K.\ Drakakis, J.\ Healy, and S.\ Rickard. ``A stochastic analysis approach in the search for Costas arrays.'' International Journal of Applied Mathematics and Engineering Sciences, Volume 2, Issue 1, January-June 2008, pp.\ 73--88.
\bibitem{DIR} K.\ Drakakis, F.\ Iorio, and S.\ Rickard. ``The enumeration of Costas arrays of order 28 and its consequences.'' Advances in Mathematics of Communications, Volume 5, Issue 1, 2011, pp.\ 69--86.
\bibitem{DIRW} K.\ Drakakis, F.\ Iorio, S.\ Rickard, and J.\ Walsh. ``Results of the enumeration of Costas arrays of order 29.'' (Advances in Mathematics of Communications.)
\bibitem{DMR} K.\ Drakakis, G.\ McGuire, and V.\ Requena. ``On the Nonlinearity of Exponential Welch Costas Functions.'' IEEE Transactions on Information Theory, Volume 56 , Issue 3, March 2010, pp.\ 1230--1238.
\bibitem{DR1} K.\ Drakakis and S.\ Rickard. ``On the construction of nearly optimal Golomb rulers by unwrapping Costas arrays.'' Contemporary Engineering Sciences, Volume 3, Issue 7, 2010, pp.\ 295--309.
\bibitem{DR2} K.\ Drakakis and S.\ Rickard. ``On the generalization of the Costas property in the continuum.'' Advances in Mathematics of Communications, Volume 2, Issue 2, 2008, pp.\ 113--130.
\bibitem{DRBCIOW} K.\ Drakakis, S.\ Rickard, J.\ Beard, R.\ Caballero, F.\ Iorio, G.\ O'Brien, and J.\ Walsh. ``Results of the enumeration of Costas arrays of order 27.'' IEEE Transactions on Information Theory, Volume 54, Issue 10, October 2008, pp.\ 4684--4687.
\bibitem{DRG} K.\ Drakakis, S.\ Rickard, and R.\ Gow. ``Interlaced Costas arrays do not exist.'' Mathematical Problems in Engineering, Volume 2008.
\bibitem{DT} D.\ Drumheller and E.\ Titlebaum. ``Cross-Correlation Properties of Algebraically Constructed Costas Arrays.'' IEEE Transactions on Aerospace and Electronic Systems, Volume 27, Issue 1, pp.\ 2--10, January 1991.
\bibitem{EGRT} P.\ Erd\"os, R.\ Graham, I.Z.\ Ruzsa, and H.\ Taylor. ``Bounds for arrays of dots with distinct slopes or lengths.'' Combinatorica, Volume 12, Issue 1, pp.\ 39--44, 1992.
\bibitem{ET} P.\ Erd\"os and P.\ Turan. ``On a problem of Sidon in additive number theory and some related problems.'' Journal of the London Mathematical Society, Volume 16, 1941, pp.\ 212--215 --- followed by Addendum (by P. Erd\"os), ibid.\ Volume 19, 1944, pp.\ 208.
\bibitem{E} T.\ Etzion. ``Combinatorial designs with Costas arrays properties.'' Discrete Mathematics, Volume 93, pp.\ 143--154, 1991.
\bibitem{E2} T. Etzion. ``Sequence folding, lattice tiling, and multidimensional coding.'' (in www.arxiv.org)
\bibitem{FL} A.\ Freedman and N.\ Levanon. ``Any two $N\times N$ Costas signals must have at least one common ambiguity sidelobe if $N>3$---A proof.'' Proceedings of the IEEE, Volume 73, Issue 10, pp. 1530--1531, October 1985.
\bibitem{G0} R.A.\ Games. ``An algebraic construction of SONAR sequences using M-sequences.'' SIAM Journal on Algebraic and Discrete Methods, Volume 8, Issue 4, pp.\ 753--761, October 1987.
\bibitem{GJ} M.R.\ Garey and D.S.\ Johnson. ``Computers and Intractability: A Guide to the Theory of NP-Completeness.'' W.H.\ Freeman \& Co., 1979.
\bibitem{G} S.\ Golomb. ``Algebraic constructions for Costas arrays.'' Journal Of Combinatorial Theory Series A, Volume 37, Issue 1, pp.\ 13--21, 1984.
\bibitem{G2} S.\ Golomb. ``The $T_4$ and $G_4$ constructions for Costas arrays.'' IEEE Transactions on Information Theory, Volume 38, Issue 4, pp.~1404--1406, July 1992.
\bibitem{G3} S.\ Golomb. ``How to number a graph.'' Graph Theory and Computing (R.C.\ Read, editor), Academic Press, New York, 1977, pp.\ 23–-37.
\bibitem{GG} S.\ Golomb and G.\ Gong. ``The status of Costas arrays.'' IEEE Transactions on Information Theory, Volume 53, Issue 11, pp.\ 4260--4265, November 2007.
\bibitem{GG2} S.\ Golomb and G.\ Gong. ``Signal design for good correlation.'' Cambridge, 2005.
\bibitem{GM} S.\ Golomb and O.\ Moreno. ``On periodicity properties of Costas arrays and a conjecture on permutation polynomials.''  IEEE Transactions on Information Theory, Volume 42, Issue 6, November 1996, pp.\ 2252--2253.
\bibitem{GT} S.\ Golomb and H.\ Taylor. ``Constructions and properties of Costas arrays'', Proceedings of the IEEE, Vol. 72, pp.\ 1143--1163, 1984.
\bibitem{HZ} J.\ Hamkins and K.\ Zeger. ``Improved bounds on maximum size binary RADAR arrays," IEEE Transactions on Information Theory, Volume 43, Issue 3, pp.\ 997--1000, May 1997.
\bibitem{HM} J.\ Holden and P.\ Moree. ``Some heuristics and results for small cycles of the discrete logarithm.'' Mathematics of Computation, Volume 75, Number 253, pp.\ 419--449, 2005.
\bibitem{IR} K.\ Ireland and M.\ Rosen. ``A Classical Introduction to Modern Number Theory.'' $2^{nd}$ edition, Springer, 1990.
\bibitem{LT} H.\ Lefmann and T.\ Thiele. ``Point sets with distinct distances.'' Combinatorica, Volume 15, Issue 3, pp.\ 379--408, 1995.
\bibitem{LPS} M.\ Levin, C.\ Pomerance, and K.\ Soundararajan. ``Fixed points for discrete logarithms.'' Proceedings of the Ninth Algorithmic Number Theory Symposium (ANTS IX), Springer Lectures Notes in Computer Science 6197, pp.\ 6--15, 2010.
\bibitem{L} B.\ Lindstr\"om. ``An inequality for $B_2$-sequences.'' Journal of Combinatorial Theory, Volume 6, 1969, pp.\ 211--212.
\bibitem{L2} B.\ Lindstr\"om. ``Finding finite $B_2$-sequences faster.'' Mathematics of Computation, Volume 67, 1998, pp.\ 1173--1178.
\bibitem{MST} S.\ Maric, I.\ Seskar, and E.\ Titlebaum. ``On Cross-Ambiguity Properties of Welch-Costas Arrays When Applied in SS/FH Multiuser Radar and Sonar Systems.'' IEEE Transactions on Aerospace and Electronic Systems, Volume 30, Issue 4, pp.\ 489--493, October 1994.
\bibitem{MGC} O.\ Moreno, S.W.\ Golomb, and C.\ Corrada. ``Extended SONAR sequences.'' IEEE Transactions on Information Theory, Volume 43, Issue 6, pp.\ 1999--2005, November 1997.
\bibitem{MGT} O. Moreno, R.A.\ Games, and H.\ Taylor. ``SONAR sequences from Costas arrays and the best known sonar sequences with up to 100 symbols.'' IEEE Transactions on Information Theory, Volume 39, Issue 6, pp.\ 1985--1987, September 1993.
\bibitem{PT} R.E.\ Peile and H.\ Taylor. ``Sets of points with pairwise distinct slopes.'' Computers and Mathematics with Applications, Volume 39, pp.\ 109--115, 2000.
\bibitem{P} B.\ Popovic. ``New constructions of Costas sequences.'' Electronics Letters, Volume 25, Issue 1, pp.\ 40--41, January 1989.
\bibitem{R} S.\ Rickard. ``Large sets of frequency hopped waveforms with nearly ideal orthogonality properties.'' Masters thesis, MIT, 1993.
\bibitem{R2} S.\ Rickard. ``Searching for Costas arrays using periodicity properties.'' IMA International Conference on Mathematics in Signal Processing at The Royal Agricultural College,'' Cirencester, UK, 2004.
\bibitem{R3} S.\ Rickard. ``Open problems in Costas arrays.'' IMA International Conference on Mathematics in Signal Processing at The Royal Agricultural College,'' Cirencester, UK, 2006.
\bibitem{RCDLW} S.\ Rickard, E.\ Connell, F.\ Duignan, B.\ Ladendorf and A.\ Wade. ``The enumeration of Costas arrays of size 26.'' Conference on Information Signals and Systems, Princeton University, USA, 2006.
\bibitem{RH} S.\ Rickard and J.\ Healy. ``Stochastic search for Costas arrays.'' Conference on Information Signals and Systems, Princeton University, USA, 2006.
\bibitem{R4} J.P.\ Robinson. ``Golomb rectangles as folded rulers.'' IEEE Transactions on Information Theory, Volume 43, Issue 1,  Janueary 1997, pp.\ 290--293.
\bibitem{R5} J.P.\ Robinson. ``Golomb rectangles.'' IEEE Transactions on Information Theory, Volume 31, Issue 6, November 1985, pp.\ 781--787.
\bibitem{R7} J.P.\ Robinson. ``Genetic Search for Golomb Arrays.'' IEEE Transactions On Information Theory, Volume 46, Issue 3, May 2000, pp.\ 1170--1173.
\bibitem{R8} J.P.\ Robinson. ``Optimum Golomb rulers.'' IEEE Transactions On Computers, Volume 28, Issue 12, December 1979, pp.\ 943--944.
\bibitem{RB} J.P.\ Robinson and A.J.\ Bernstein. ``A class of binary recurrent codes with limited error propagation.'' IEEE Transactions on  Information Theory, Volume 13, Issue 1, 1967, pp.\ 106--113.
\bibitem{R6} I.Z.\ Ruzsa. ``Solving a linear equation in a set of integers I.'' Acta Arithmetica, Volume LXV, Issue 3, 1993, pp.\ 259--282.
\bibitem{S3} J.B.\ Shearer's webpage on Golomb rulers:\newline http://www.research.ibm.com/people/s/shearer/grtab.html.
\bibitem{S2} S.\ Sidon. ``Ein Satz \"uber trigonometrische Polynome und seine Anwendungen in der Theorie der Fourier-Reihen.'' Mathematische Annalen, Volume 106, 1932, pp.\ 536--539 (in German).
\bibitem{S1} J.\ Singer. ``A theorem in finite projective geometry and some applications to number theory.'' Transactions of the American Mathematical Society, Volume 43, 1938, pp.\ 377--385.
\bibitem{S} R.P.\ Stanley. ``Enumerative combinatorics, Volume I.'' Cambridge, 1997.
\bibitem{SVM} J.\ Silverman, V.\ Vickers, and J.\ Mooney. ``On the Number of Costas arrays as a function of array size.'' Proceedings of the IEEE, Volume 76, Issue 7, pp.\ 851--853, July 1988.
\bibitem{TM} E.\ Titlebaum and S.\ Maric. ``Multiuser SONAR properties for Costas array frequency hop coded signals.'' In Proceedings of IEEE ICASSP 1990, pp.\ 2727--2730, 1990.
\bibitem{TRD} K.\ Taylor, S.\ Rickard, and K.\ Drakakis. ``Costas Arrays: Survey, Standardization and MATLAB Toolbox.'' ACM Transactions on Mathematical Software, Volume 37, Issue 4, Dec 2010.
\bibitem{W} R.A.\ Wilson. ``The finite simple groups.'' Graduate Texts in Mathematics 251, Springer-Verlag, 2009.
\bibitem{YOH} Z.\ Ye, J.\ Ouyang, and R.\ He. ``The discovery of $53\times 53$ Costas arrays.'' Natural Science Journal of Xiangtan University, Volume 17, Issue 4, Dec 1995, pp.\ 120--122 (in Chinese).
\bibitem{Z} W.P.\ Zhang. ``On a problem of Brizolis.'' Pure and Applied Mathematics, Volume 11 (suppl.), pp.\ 1--3, 1995 (in Chinese).
\bibitem{Z2} W.P.\ Zhang. ``A note on arrays of dots with distinct slopes.'' Combinatorica, Volume 13, Issue 1, pp.\ 127--128, 1993.
\bibitem{ZT} Z.\ Zhang and C.\ Tu. ``New bounds for the sizes of RADAR arrays.'' IEEE Transactions on Information Theory, Volume 40, Issue 5, pp.\ 1672--1678, September 1994.
\bibitem{cawp} http://www.costasarrays.org
\bibitem{cawpb} http://jameskbeard.com/jameskbeard/
\bibitem{ogr} http://www.distributed.net/OGR/en
\end{thebibliography}
\end{document}